\documentclass{amsart}
\usepackage{amssymb,amscd}
\usepackage[all]{xy}
\usepackage[stdtext,vcentermath]{myyoungtab}

\numberwithin{equation}{section}

\newtheorem{prop}{Proposition}
\newtheorem{theorem}[prop]{Theorem}
\newtheorem{corollary}[prop]{Corollary}
\newtheorem{lemma}[prop]{Lemma}
\newtheorem{conjecture}[prop]{Conjecture}

\theoremstyle{definition}
\newtheorem{definition}[prop]{Definition}
\newtheorem{example}[prop]{Example}
\newtheorem{remark}[prop]{Remark}

\numberwithin{prop}{section}

\title{$X=M$ for symmetric powers}

\author{Anne Schilling}
\address{Department of Mathematics \\
University of California \\
One Shields Ave. \\
Davis, CA 95616-8633 U.S.A. }
\email{anne@math.ucdavis.edu}

\author{Mark Shimozono}
\address{Department of Mathematics \\
Virginia Tech \\
Blacksburg, VA 24061-0123 USA} \email{mshimo@math.vt.edu}

\thanks{\textit{Date:} December 2004}
\thanks{AS was supported in part by NSF grant DMS-0200774.}
\thanks{MS was supported in part by NSF grant DMS-0100918.}

\newcommand{\phib}{\overline{\phi}}
\newcommand{\phit}{\widetilde{\phi}}

\newcommand{\flip}{\theta}

\newcommand{\lh}{\mathrm{lh}}
\newcommand{\rh}{\mathrm{rh}}
\newcommand{\ldh}{\lh^\vee}

\newcommand{\Lh}{\lh(L)}
\newcommand{\Lr}{\rh(L)}
\newcommand{\Ldh}{\ldh(L)}

\newcommand{\ls}{\mathrm{ls}}
\newcommand{\rs}{\mathrm{rs}}
\newcommand{\lds}{\ls^\vee}

\newcommand{\Ls}{\ls(L)}
\newcommand{\Lrs}{\rs(L)}
\newcommand{\Lds}{\lds(L)}

\newcommand{\llb}{\mathrm{lb}}
\newcommand{\rb}{\mathrm{rb}}

\newcommand{\Lb}{\llb(L)}

\newcommand{\db}{\overline{\delta}}
\newcommand{\ddb}{\db^\vee}
\newcommand{\dt}{\widetilde{\delta}}

\newcommand{\rcls}{\overline{j}}
\newcommand{\rcrs}{\widetilde{j}}
\newcommand{\rclds}{\rcls^\vee}

\newcommand{\rclb}{\overline{i}}
\newcommand{\rcrb}{\widetilde{i}}

\newcommand{\0}{\circ}

\newcommand{\aut}{\sigma}

\newcommand{\CC}{\mathcal{C}}
\newcommand{\CCA}{\CC^A}
\newcommand{\CCAD}{\CC^{A\vee}}
\newcommand{\Cf}{\mathcal{C}^{\mathrm{fin}}}

\newcommand{\level}{\mathrm{level}}
\newcommand{\ve}{\varepsilon}
\newcommand{\vp}{\varphi}

\newcommand{\na}{\natural}

\newcommand{\wt}{\mathrm{wt}}
\newcommand{\Xt}{\overleftarrow{X}}

\newcommand{\es}{\varnothing}
\newcommand{\vn}{\es}
\newcommand{\cc}{cc}
\newcommand{\Conf}{\mathrm{C}}
\newcommand{\lb}{\bar{\ell}}
\newcommand{\sbar}{\bar{s}}

\newcommand{\nh}{\widehat{\nu}}
\newcommand{\Jh}{\widehat{J}}
\newcommand{\mh}{\widehat{m}}

\newcommand{\RC}{\mathrm{RC}}
\newcommand{\RCt}{\widetilde{\RC}}
\newcommand{\RCtv}{\RCt^v}
\newcommand{\qbin}[2]{\genfrac{[}{]}{0pt}{}{#1}{#2}}

\newcommand{\rkt}{\widetilde{\rk}}
\newcommand{\rk}{\mathrm{rk}}
\newcommand{\vrk}{\rk^v}
\newcommand{\rkh}{\widehat{\rk}}

\newcommand{\Lp}{\lh^{-1}(L)}
\newcommand{\Lhat}{\widehat{L}}

\newcommand{\Dt}{\overleftarrow{D}}

\newcommand{\emb}{\Psi}
\newcommand{\bij}{\iota}
\newcommand{\mult}{\gamma}

\newcommand{\eh}{\hat{e}}
\newcommand{\fh}{\hat{f}}

\newcommand{\Pv}{P^v}
\newcommand{\RCv}{\RC^v}
\newcommand{\Rh}{\hat{R}}
\newcommand{\Vh}{\widehat{V}}

\newcommand{\vrs}{\widehat{\rs}}

\newcommand{\vls}{\widehat{\ls}}
\newcommand{\vlh}{\widehat{\lh}}
\newcommand{\vdb}{\widehat{\delta}}
\newcommand{\vrcls}{\widehat{j}}
\newcommand{\vphib}{\phib^v}

\newcommand{\gehb}{\overline{\geh}}
\newcommand{\geh}{\mathfrak{g}}
\newcommand{\J}{\bar{I}}

\newcommand{\Wb}{\overline{W}}

\newcommand{\Pb}{\overline{P}}
\newcommand{\Pcl}{P'}

\newcommand{\inner}[2]{\langle #1\,,\,#2\rangle}

\newcommand{\HH}{\mathcal{H}}

\newcommand{\Lab}{\overline{\Lambda}}
\newcommand{\La}{\Lambda}
\newcommand{\la}{\lambda}
\newcommand{\lm}{\lambda^-}
\newcommand{\lp}{\lambda^+}

\newcommand{\Image}{\mathrm{Im}}

\newcommand{\fst}{\mathrm{st}}
\newcommand{\nth}{\mathrm{th}}

\newcommand{\Z}{\mathbb{Z}}

\newcommand{\wsig}{\tau}

\newcommand{\bh}{\widehat{b}}
\newcommand{\xh}{\widehat{x}}
\newcommand{\bhp}{\bh'}

\newcommand{\ba}{\bar{1}}
\newcommand{\bb}{\bar{2}}
\newcommand{\bc}{\bar{3}}

\newcommand{\be}{\bar{5}}

\newcommand{\xx}{{\hphantom{x}}}

\newcommand{\ov}[1]{\overline{#1}}

\begin{document}

\begin{abstract}
The $X=M$ conjecture of Hatayama et al. asserts the equality between
the one-dimensional configuration sum $X$ expressed as the generating
function of crystal paths with energy statistics and the fermionic formula $M$
for all affine Kac--Moody algebra. In this paper we prove the
$X=M$ conjecture for tensor products of Kirillov--Reshetikhin crystals
$B^{1,s}$ associated to symmetric powers for all nonexceptional affine algebras.
\end{abstract}

\maketitle


\section{Introduction}
In two extraordinary papers, Hatayama et al.~\cite{HKOTT:2001,HKOTY:1999}
recently conjectured the equality between the one-dimensional configuration
sum $X$ and the fermionic formula $M$ for all affine Kac--Moody algebras.
The one-dimensional configuration sum $X$ originates from the corner-transfer-matrix
method~\cite{Baxter:1982} used to solve exactly solvable lattice models in statistical
mechanics. It is the generating function of highest weight crystal paths
graded by the energy statistic. The fermionic formula $M$ comes from the Bethe
Ansatz~\cite{Bethe:1931} and exhibits the quasiparticle structure of the underlying
model. In combinatorial terms, it can be written as the generating function of
rigged configurations.

The one-dimensional configuration sum depends on the underlying tensor
product of crystals. In~\cite{HKOTT:2001,HKOTY:1999}, the $X=M$ conjecture was
formulated for tensor products of Kirillov--Reshetikhin (KR) crystals $B^{r,s}$.
Kirillov--Reshetikhin crystals are crystals for finite-dimensional irreducible
modules over quantum affine algebras. The irreducible finite-dimensional
$U'_q(\geh)$-modules were classified by Chari and Pressley~\cite{CP:1995,CP:1998}
in terms of Drinfeld polynomials.  The Kirillov--Reshetikhin modules $W^{r,s}$,
labeled by a Dynkin node $r$ of the underlying classical algebra and a positive
integer $s$, form a special class of these
finite-dimensional modules.  They naturally correspond to the weight
$s\Lambda_r$, where $\Lambda_r$ is the $r$-th fundamental weight of
$\geh$. It was conjectured in~\cite{HKOTT:2001,HKOTY:1999}, that there exists a crystal
$B^{r,s}$ for each $W^{r,s}$. In general, the existence of $B^{r,s}$ is still an
open question. For type $A_n^{(1)}$ the crystal $B^{r,s}$ is known to
exist~\cite{KKMMNN:1992} and its combinatorial structure has been studied~\cite{Sh:2002}.
The crystals $B^{1,s}$ for nonexceptional types, which are relevant for this paper,
are also known to exist and their combinatorics has been worked
out~\cite{KKM:1994,KKMMNN:1992}.

The purpose of this paper is to establish the $X=M$ conjecture for
tensor products of KR crystals of the form $B^{1,s}$ for
nonexceptional affine algebras. This extends~\cite{OSS:2002a},
where $X=M$ is proved for tensor powers of $B^{1,1}$,
and~\cite{KR:1988,KSS:2002}, where $X=M$ is proved for type $A_n^{(1)}$.

Our method to prove $X=M$ for symmetric powers combines various
previous results and techniques. $X=M$ is first proved for $\geh$
such that $\gehb$ is simply-laced (see Corollary~\ref{cor:X=M AD}).
This is accomplished by exhibiting a grade-preserving bijection from
$U'_q(\gehb)$-highest weight vectors (paths) to rigged
configurations (RCs). This was already proved for the root system
$A_n^{(1)}$ \cite{KSS:2002}. For type $D_n^{(1)}$ we exhibit such a
path-RC bijection. The proof essentially reduces to the previously
known $s=1$ case \cite{OSS:2002a} using the ``splitting" maps
$B^{1,s}\rightarrow B^{1,s-1}\otimes B^{1,1}$ which are
$U_q(\gehb)$-equivariant grade-preserving embeddings.

To prove that the bijection preserves the grading, we consider an
involution denoted $*$ on crystal graphs that combines
contragredient duality with the action of the longest element $w_0$
of the Weyl group of $\gehb$. This duality on the crystal graph,
corresponds under the path-RC bijection to the involution on RCs
given by complementing the quantum numbers with respect to the
vacancy numbers.

We then reduce to the case that $\gehb$ is simply-laced. This is
achieved using the embedding of an affine algebra $\geh$ into one
(call it $\geh_Y$) whose canonical simple Lie subalgebra is
simply-laced. On the $X$ side we use the virtual crystal
construction developed in \cite{OSS:2003a,OSS:2003b}. It is shown in
\cite{OSS:2003b} that the KR $U'_q(\geh)$-crystals $B^{1,s}$ embed
into tensor products of KR $U'_q(\geh_Y)$-crystals such that the
grading is respected. One may define the $VX$ (``virtual X") formula
in terms of the image of this embedding and show that $X=VX$ (see
section~\ref{ss:VX=V}). This is proved for tensor products of
crystals $B^{1,s}$ in~\cite{OSS:2003b}. On the $M$ side, it is
observed in \cite{OSS:2003b} that the RCs giving the fermionic
formula $M$ for type $\geh$, embed into the set of RCs giving a
fermionic formula for type $\geh_Y$. Let us denote by $VM$
(``virtual $M$") the generating function over the image of this
embedding of fermionic formulas. It is shown in \cite{OSS:2003b}
that $M=VM$. It then suffices to prove $VX=VM$. That is, one must
show that the path-to-RC bijection that has already been established
for the simply-laced cases, restricts to a bijection between the
subsets of objects in the formulas $VX$ and $VM$. This is shown in
Theorem~\ref{th:virtualbij} and as a corollary proves $X=M$ for
nonsimply-laced algebras, as stated in
Corollary~\ref{cor:X=VX=VM=M}.

In section~\ref{s:X} we review the crystal theory, the definition of
the one-dimensional configuration sum $X$, contragredient duality
and the $*$ involution. Virtual crystals are reviewed in
section~\ref{s:virtual}. Right and left splitting of crystals are
discussed in sections~\ref{s:right} and~\ref{s:left}, respectively.
Rigged configurations and the analogs of the splitting maps are
subject of section~\ref{sec:RC}. The fermionic formulas $M$ and
their virtual counterparts $VM$ are stated in
section~\ref{s:fermionic}. The $X=M$ conjecture for types
$A_n^{(1)}$ and $D_n^{(1)}$ is proven in section~\ref{sec:bij} by
establishing a statistics preserving bijection. Finally, in
section~\ref{sec:virtualbij} the equality $X=VX=VM=M$ is established
for nonsimply-laced types.

\section{Formula $X$} \label{s:X}

\subsection{Affine algebras}
Let $\geh \supset \geh' \supset \gehb$ be a nonexceptional affine
Kac-Moody algebra, its derived subalgebra and canonical simple Lie
subalgebra \cite{Kac:1985}. Denote the corresponding quantized
universal enveloping algebras by $U_q(\geh)\supset
U_q'(\geh)\supset U_q(\gehb)$ \cite{Ka:95}. Let $I=\J\cup \{0\}$
(resp. $\J$) be the vertex set of the Dynkin diagram of $\geh$
(resp. $\gehb$). For $i\in I$, let $\alpha_i$, $h_i$, $\La_i$ be
the simple roots, simple coroots, and fundamental weights of
$\geh$. Let $\{\Lab_i\mid i\in \J\}$ be the fundamental weights of
$\gehb$. Let $(a_0,a_1,\dotsc,a_n)$ be the smallest tuple of
positive integers giving a dependency relation on the columns of
the Cartan matrix of $\geh$. Write $a_i^\vee$ for the
corresponding integers for the Langlands dual Lie algebra, the one
whose Cartan matrix is the transpose of that of $\geh$. Let
$c=\sum_{i\in I} a_i^\vee h_i$ be the canonical central element
and $\delta=\sum_{i\in I} a_i \alpha_i$ the generator of null
roots. Let $Q,Q^\vee,P$ be the root, coroot, and weight lattices
of $\geh$. Let $\inner{\cdot}{\cdot}:Q^\vee\otimes P \rightarrow
\Z$ be the pairing such that $\inner{h_i}{\La_j}=\delta_{ij}$. Let
$P\rightarrow \Pcl\rightarrow\Pb$ be the natural surjections of
weight lattices of $\geh\supset \geh'\supset \gehb$. Let
$\Pb^+\subset \Pb$ be the dominant weights for $\gehb$. Let $W$
and $\Wb$ be the Weyl groups of $\geh$ and $\gehb$ respectively.

\subsection{Crystal graphs}
Let $M$ be a finite-dimensional $U'_q(\geh)$-module. Such modules
are not highest weight modules (except for the zero module) and
therefore need not have a crystal base. Suppose $M$ has a crystal
base $B$. This is a special basis of $M$; it possesses the
structure of a colored directed graph called the \textbf{crystal graph}. By
abuse of notation the vertex set of the crystal graph is also
denoted $B$. The edges of the crystal graph are colored by the set
$I$. It has the following properties (that of a regular
$\Pb$-weighted $I$-crystal):
\begin{enumerate}
\item Fix an $i\in I$. If all edges are removed except those colored $i$, the
connected components are finite directed linear paths called the
\textbf{$i$-strings} of $B$. Given $b\in B$, define $f_i(b)$
(resp. $e_i(b)$) to be the vertex following (resp. preceding) $b$
in its $i$-string; if there is no such vertex, declare the result
to be the special symbol $\emptyset$. Define $\vp_i(b)$ (resp.
$\ve_i(b)$) to be the number of arrows from $b$ to the end (resp.
beginning) of its $i$-string.
\item There is a function $\wt:B\rightarrow \Pb$ such that
\begin{equation*}
\begin{split}
\wt(f_i(b))&=\wt(b)-\alpha_i \\
\vp_i(b)-\ve_i(b) &= \inner{h_i}{\wt(b)}.
\end{split}
\end{equation*}
\end{enumerate}

A \textbf{morphism} $g:B\rightarrow B'$ of $\Pb$-weighted
$I$-crystals is a map $g:B\cup\{\emptyset\}\rightarrow
B'\cup\{\emptyset\}$ such that $g(\emptyset)=\emptyset$ and for
any $b\in B$ and $i\in I$, $g(f_i(b))=f_i(g(b))$ and
$g(e_i(b))=e_i(g(b))$. An isomorphism of crystals is a morphism of
crystals which is a bijection whose inverse bijection is also a
morphism of crystals.

If $B_i$ is the crystal base of the $U'_q(\geh)$-module $M_i$ for
$i=1,2$ then the \textbf{tensor product} $M_2 \otimes M_1$ is a
$U'_q(\geh)$-module with crystal base denoted $B_2\otimes B_1$.
Its vertex set is just the cartesian product $B_2 \times B_1$. Its
edges are given in terms of those of $B_1$ and $B_2$ as follows.

\begin{remark}
We use the opposite of Kashiwara's tensor product convention.
\end{remark}
One has $\wt(b_2\otimes b_1)=\wt(b_2)+\wt(b_1)$ and
\begin{equation*}
\begin{split}
  f_i(b_2\otimes b_1) &= \begin{cases}
  f_i(b_2)\otimes b_1 &\text{if $\ve_i(b_2)\ge\vp_i(b_1)$} \\
  b_2 \otimes f_i(b_1) &\text{otherwise,}
\end{cases} \\
  e_i(b_2\otimes b_1) &= \begin{cases}
  e_i(b_2)\otimes b_1 &\text{if $\ve_i(b_2)>\vp_i(b_1)$} \\
  b_2 \otimes e_i(b_1) &\text{otherwise,}
\end{cases}
\end{split}
\end{equation*}
where the result is declared to be $\emptyset$ if either of its
tensor factors are.

The tensor product construction is associative up to isomorphism.

Define $\vp,\ve:B\rightarrow\Pcl$ by
\begin{equation*}
  \vp(b) = \sum_{i\in I} \vp_i(b) \La_i \qquad\qquad \ve(b) = \sum_{i\in I} \ve_i(b) \La_i.
\end{equation*}

Every irreducible integrable finite-dimensional $U_q(\gehb)$-module
is a highest weight module with some highest weight $\la\in\Pb^+$;
denote its crystal graph by $B(\la)$. It is a $\Pb$-weighted
$\J$-crystal with a unique classical highest weight vector.

A \textbf{classical component} of the crystal graph $B$ of a
$U'_q(\geh)$-module is a connected component of the graph obtained
by removing all $0$-arrows from $B$. The vertex $b\in B$ is a
\textbf{classical highest weight vector} if $\ve_i(b)=0$ for all
$i\in \J$. Each classical component of a $U'_q(\geh)$-module has a
unique classical highest weight vector.

\subsection{Finite crystals} \label{ss:fin}
Let $\Cf$ be the category of \textbf{finite} crystals as defined
in \cite{HKKOT:2000}. Every $B\in \Cf$ has the following
properties.
\begin{enumerate}
\item $B$ is the crystal base of an irreducible $U'_q(\geh)$-module
and is therefore connected.
\item There is a weight $\la\in \Pb^+$ such that there
is a unique $u(B)\in B$ with $\wt(u(B))=\la$
and for all $b\in B$, $\wt(b)$ is in the convex hull of $\Wb \la$.
\end{enumerate}
$\Cf$ is a tensor category \cite{HKKOT:2000}. If $B,B'\in \Cf$
then $B \otimes B'\in \Cf$ is connected and $u(B\otimes
B')=u(B)\otimes u(B')$. Due to the existence of the universal
$R$-matrix for $U'_q(\gehb)$ it follows from \cite{KMN:1992} that:
\begin{enumerate}
\item There is a unique $U'_q(\geh)$-crystal isomorphism
$R_{B,B'}:B \otimes B'\rightarrow B'\otimes B$ called the
\textbf{combinatorial $R$-matrix}.
\item There is a unique
function (the \textbf{local coenergy}) $H=H_{B,B'}:B \otimes
B'\rightarrow \Z_{\ge0}$ that is constant on classical components,
zero on $u(B\otimes B')$, and is such that if $R_{B,B'}(b \otimes
b')=c'\otimes c$ then
\begin{equation} \label{eq:H}
  H(e_0(b\otimes b')) = H(b\otimes b') +
\begin{cases}
  1 &
\begin{aligned}
\text{if } e_0(b\otimes b')&=e_0(b)\otimes b' \text{ and } \\
e_0(c'\otimes c)&=e_0(c')\otimes c
\end{aligned}  \\
-1 &
\begin{aligned}
\text{if } e_0(b\otimes b')&=b \otimes e_0(b') \text{ and } \\
e_0(c'\otimes c)&=c'\otimes e_0(c)
\end{aligned}  \\
0 &
\text{otherwise.}
\end{cases}
\end{equation}
\end{enumerate}
The combinatorial $R$-matrices satisfy
\begin{equation*}
\begin{split}
  R_{B,B} &= 1_{B\otimes B} \\
  R_{B_1,B_2} \circ R_{B_2,B_1} &= 1_{B_2\otimes B_1}
\end{split}
\end{equation*}
and the Yang Baxter equation, the equality of isomorphisms
$B_3\otimes B_2 \otimes B_1\rightarrow B_1 \otimes B_2 \otimes
B_3$ given by
\begin{equation} \label{eq:RYB}
\begin{split}
 &(1_{B_1} \otimes R_{B_3,B_2})\circ (R_{B_3,B_1}\otimes 1_{B_2})
\circ (1_{B_3} \otimes R_{B_2,B_1}) \\
=\, &(R_{B_2,B_1} \otimes 1_{B_3}) \circ (1_{B_2} \otimes
R_{B_3,B_1}) \circ (R_{B_3,B_2} \otimes 1_{B_1}).
\end{split}
\end{equation}
We shall abuse notation and write $R_j$ (resp. $H_j$) to denote
the application of an appropriate combinatorial $R$-matrix (resp. local
coenergy function) on the $(j+1)$-th and $j$-th tensor factors
\textit{from the right}. Then \eqref{eq:RYB} reads
$R_1 R_2 R_1=R_2 R_1 R_2$. One has the following identities on a three-fold
tensor product:
\begin{equation*}
\begin{split}
  H_2 + H_1 R_2 &= H_2 R_1 + H_1 R_2 R_1 \\
  H_1 + H_2 R_1 &= H_1 R_2 + H_2 R_1 R_2.
\end{split}
\end{equation*}

\begin{prop} \label{pp:Hbraided} \cite{OSS:2003a}
Let $B=B_L\otimes\dotsm\otimes B_1$ and $B'=B_M'\otimes\dotsm \otimes B_1'$.
\begin{enumerate}
\item $R_{B,B'}$ is equal to any composition of $R$-matrices of the
form $R_{B_i,B_j'}$ which shuffle the $B_i$ to the right, past the
$B_j'$.
\item For $b\otimes b'\in B\otimes B'$, the value of $H_{B,B'}$ is the sum of the values
$H_{B_i \otimes B_j'}$ evaluated at the pairs of elements in $B_i
\otimes B_j'$ that must be switched by an $R$-matrix $R_{B_i
,B_j'}$ in the computation of $R_{B,B'}(b\otimes b')$.
\end{enumerate}
\end{prop}

\subsection{Categories $\CC$ and $\CCA$ of KR crystals}
\label{ss:tensorcats} We work with two categories of crystals. Let
$\geh$ be of nonexceptional affine type. The KR modules $W^{(1)}_s$
and their crystal bases $B^s:=B^{1,s}$ were constructed in
\cite{KKM:1994}. See also \cite{OSS:2003b} for an explicit
description of $B^s$. Let $\CC$ be the category of tensor products
of KR crystals of the form $B^s$. One has that $\CC\subset \Cf$.

The crystal $B^s$ has the $U_q(\gehb)$-decomposition
\begin{equation} \label{eq:ClDecomp}
  B^s \cong \begin{cases}
B(s\Lab_1) & \text{for $A_n^{(1)}$, $B_n^{(1)}$, $D_n^{(1)}$, $A_{2n-1}^{(2)}$} \\
\displaystyle{\bigoplus_{r=0}^s B((s-r)\Lab_1)}
&\text{for $A_{2n}^{(2)},D_{n+1}^{(2)}$} \\
\displaystyle{\bigoplus_{r=0}^{\lfloor \frac{s}{2}\rfloor}
B((s-2r)\Lab_1)} & \text{for $C_n^{(1)}$, $A_{2n}^{(2)\dagger}$.}
\end{cases}
\end{equation}
In particular $u(B^s)$ is the unique vector of weight $s\Lab_1$ in
$B^s$.

Let $\CCA$ be the category of all tensor products of KR crystals
$B^{r,s}$ in type $A_n^{(1)}$. Here $B^{r,s} \cong B(s\Lab_r)$. So
$u(B^{r,s})$ is the unique vector in $B^{r,s}$ of weight $s\Lab_r$.
$B^{r,s}$ consists of the semistandard Young tableaux of shape given
by an $r\times s$ rectangle, with entries in the set
$\{1,2,\dotsc,n+1\}$ \cite{KN:1994}. The structure of $B^{r,s}$ as
an affine crystal was given explicitly in \cite{Sh:2002}.

We fix some notation for $B\in \CC$ or $B\in\CCA$. Let $\HH=\J\times
\Z_{>0}$ where recall that $\J=\{1,2,\dotsc,n\}$ is the set of Dynkin nodes for
$\gehb$. The multiplicity array of $B$ is the array
$L=(L_i^{(a)}\mid(a,i)\in\HH)$ such that $L_i^{(a)}$ is the number
of times $B^{a,i}$ occurs as a tensor factor in $B$ for all
$(a,i)\in \HH$. Up to reordering of tensor factors
$B=\bigotimes_{(a,i)\in\HH} (B^{a,i})^{\otimes L_i^{(a)}}$.

\subsection{Intrinsic coenergy} \label{ss:D}
For $B\in\Cf$, say that $D:B\rightarrow \Z$ is an
\textbf{intrinsic coenergy} function for $B$ if $D(u(B))=0$, $D$
is constant on $U_q(\gehb)$-components, and
\begin{equation*}
   D(e_0(b)) - D(b) \le 1\qquad\text{for all $b\in B$.}
\end{equation*}
A \textbf{graded crystal} is a pair $(B,D)$ where $B\in \Cf$ and
$D$ is an intrinsic coenergy function on $B$.

We shall give each $B\in\CC$ a particular graded crystal
structure.

For $B \in \Cf$ define
\begin{equation*}
  \level(B) = \min \{ \inner{c}{\vp(b)} \mid b\in B \}.
\end{equation*}
One may verify that there is a unique element $b^\na\in B^s$ such that
\begin{equation*}
\vp(b^\na)=\level(B^s)\La_0.
\end{equation*}
Define the intrinsic coenergy function $D_{B^s}:B^s\rightarrow\Z$
by
\begin{equation*}
  D_{B^s}(b) = H_{B^s,B^s}(b \otimes b^\na)-H_{B^s,B^s}(u(B^s)\otimes b^\na).
\end{equation*}

\begin{example} \label{ex:DKR} $D_{B^s}$ has value $r$ on the
$r$-th summand in \eqref{eq:ClDecomp}.
\end{example}

\begin{prop} \label{pp:tensor} \cite{OSS:2003a}
Graded crystals form a tensor category as follows. If $(B_j,D_j)$
is a graded crystal for $1\le j\le L$, then their tensor product
$B=B_L\otimes\dotsm\otimes B_1$ is a graded crystal with
\begin{equation} \label{eq:DNY}
  D_B = \sum_{1\le i<j\le L} H_i R_{i+1} R_{i+2}\dotsm R_{j-1}
   + \sum_{j=1}^L D_{B_j} R_1 R_2 \dotsm R_{j-1}
\end{equation}
where $D_{B_j}$ acts on the rightmost tensor factor.
\end{prop}

\subsection{$X$ formula}
Let $(B,D)$ be a graded crystal. For $\la\in\Pb^+$ let $P(B,\la)$
be the set of classical highest weight vectors in $B$ of weight
$\la$. Define the one-dimensional sum
\begin{equation} \label{eq:X}
  X_{B,\la}(q) = \sum_{b\in P(B,\la)} q^{D_B(b)/a_0}.
\end{equation}
Recall that $a_0=1$ unless $\geh=A_{2n}^{(2)}$ in which case
$a_0=2$.

\subsection{Contragredient duality}
\label{ss:dual} Given a $U'_q(\geh)$-module $M$ with crystal base
$B$, the contragredient dual module $M^\vee$ has a crystal base
$B^\vee=\{b^\vee\mid b\in B \}$ such that
\begin{equation*}
\begin{split}
  \wt(b^\vee) &= - \wt(b) \\
  f_i(b^\vee) &= e_i(b)^\vee
\end{split}
\end{equation*}
for $i\in I$ and $b\in B$ such that $e_i(b)\not=\emptyset$.
\begin{prop} \label{pp:dualtensor}
\begin{equation*}
  (B_2\otimes B_1)^\vee \cong B_1^\vee \otimes B_2^\vee.
\end{equation*}
\end{prop}

\begin{example} \label{ex:dualA} Assume type $A_n^{(1)}$. We have
\begin{equation} \label{eq:dualKR}
  B^{r,s\vee} \cong B^{n+1-r,s}.
\end{equation}
The composite map
\begin{equation*}
  B^{r,s} \overset{\vee}{\longrightarrow} B^{r,s\vee} \cong
  B^{n+1-r,s}
\end{equation*}
is given explicitly as follows. Let $b\in B^{r,s}$. Replace each
column of $b$, viewed as a subset of $\{1,2,\dotsc,n+1\}$ of size
$r$, by the column of size $n+1-r$ given by its complement. Then
reverse the order of the columns. For $n=5$, $r=2$, and $s=3$, a
tableau $b\in B^{r,s}$ and its image in $B^{n+1-r,s}$ are given
below:
\begin{equation*}
  \young(112,346)\mapsto \young(122,334,455,566).
\end{equation*}
\end{example}

\begin{example}\label{ex:dualArow}
By definition $B^{1,1\vee}$ is defined by replacing each element of
$b\in B^{1,1}$ by an element $b^\vee$ and reversing arrows.
$B^{1,s\vee}$ can be realized by the weakly increasing words of
length $s$ in the alphabet $\{(n+1)^\vee<\dotsm<2^\vee<1^\vee\}$.
The arrow-reversing map from $B^s$ to $B^{s\vee}$ is given by taking
a word of length $s$, replacing each symbol $i$ with $i^\vee$, and
reversing.
\end{example}

\subsection{Dynkin automorphisms}
\label{ss:dynkin} Let $\sigma$ be an automorphism of the Dynkin
diagram of $\geh$. This induces isometries $\sigma:P\rightarrow P$
and $\sigma:\Pb\rightarrow\Pb$ given by
$\sigma(\La_i)=\La_{\sigma(i)}$ for $i\in I$, $\sigma(\delta)=\delta$, and
$\sigma(\Lab_i)=\Lab_{\sigma(i)}$ for $i\in \J$.

If $M$ is a $U'_q(\geh)$-module with crystal base $B$, then by
carrying out the construction of $M$ but with $i$ replaced
everywhere by $\sigma(i)$, there is a $U'_q(\geh)$-module
$M^\sigma$ with crystal base $B^\sigma$ and a bijection
$\sigma:B\rightarrow B^\sigma$ such that
\begin{equation*}
\begin{split}
  \wt(\sigma(b)) &= \sigma(\wt(b)) \\
  \sigma(e_i(b)) &= e_{\sigma(i)}(b) \\
  \sigma(f_i(b)) &= f_{\sigma(i)}(b)
\end{split}
\end{equation*}
for all $b\in B$ and $i\in I$.

In particular, if the appropriate KR modules have been constructed
then
\begin{equation*}
  (B^{r,s})^\sigma=B^{\sigma(r),s}.
\end{equation*}

\subsection{The Dynkin involution $\wsig$}
\label{ssec:sigma} We fix a canonical Dynkin automorphism $\wsig$ of
the affine Dynkin diagram in the following manner. There is a
length-preserving involution on $\Wb$ given by conjugation by the
longest element $w_0\in\Wb$. Restricting this involution to elements
of length one, one obtains an involution $\wsig$ on the set of
simple reflections $\{s_i\mid i\in \J\}$ of $\Wb$. For simplicity of
notation this can be written as an involution on the index set $\J$.
This gives an automorphism of the Dynkin diagram of $\wsig$. Call
the resulting Dynkin automorphism $\wsig$.

Explicitly, $\wsig$ is the identity except when $\gehb=A_{n-1}$
where $\wsig$ exchanges $i$ and $n-i$, and $\gehb=D_n$ with $n$ odd,
where $\wsig$ exchanges $n-1$ and $n$ and fixes all other Dynkin
nodes.

$\wsig$ may be extended to the Dynkin diagram of $\geh$ by fixing
the $0$ node. It satisfies $w_0 s_i w_0 = s_{\wsig(i)}$ for all
$i\in I$.

The automorphism $\wsig$ induces the following action on the weight
lattice $P$:
\begin{equation*}
  \wsig(\La_i) = \La_{\wsig(i)}\qquad\text{for $i\in I$.} \\
\end{equation*}
One may show that this is equivalent to
\begin{equation*}
  \wsig(\La) = - w_0 \,\La \qquad\text{for $\La\in P$.}
\end{equation*}
In particular
\begin{equation} \label{eq:sigalpha}
  \wsig(\alpha_i) = \alpha_{\wsig(i)} = - w_0\, \alpha_i.
\end{equation}

\subsection{The $*$ involution}
Let $M$ be a $U'_q(\geh)$-module with crystal base $B$. With $\wsig$
as above, define the module
\begin{equation*}
  M^* = M^{\wsig\vee}.
\end{equation*}
It has crystal base $B^*$, with elements $b^*$ for $b\in B$ such
that
\begin{equation} \label{eq:wtstar}
\wt(b^*)= w_0 \wt(b)
\end{equation}
and
\begin{equation} \label{eq:defstar}
\begin{split}
  e_i(b^*) &= f_{\wsig(i)}(b)^* \\
  f_i(b^*) &= e_{\wsig(i)}(b)^*
\end{split}
\end{equation}
for all $i\in I$.

\begin{remark} \label{rem:starcomp}
By \eqref{eq:defstar} for $i\in\J$ it follows that the map $*$ sends
classical components of $B$ to classical components of $B^*$, which
by \eqref{eq:wtstar} must have the same classical highest weight.
\end{remark}

\begin{prop} \label{pp:star} $(B_1\otimes B_2)^*\cong B_2^*\otimes B_1^*$
with $(b_1\otimes b_2)^* \mapsto b_2^* \otimes b_1^*$.
\end{prop}

\begin{conjecture} \label{conj:crystalstar} Let $B\in
\Cf$. Then there is a unique involution $*:B\rightarrow B$ such
that \eqref{eq:wtstar} and \eqref{eq:defstar} hold.
\end{conjecture}

Uniqueness follows from the connectedness of $B$ and the fact that
$u(B)$ is the unique vector in $B$ of its weight.

\begin{remark} \label{rem:Rstar}
The crystals satisfying Conjecture \ref{conj:crystalstar} form a
tensor category. Given involutions $*$ on $B_1$ and $B_2$
satisfying Conjecture \ref{conj:crystalstar}, define $*$ on $B_1
\otimes B_2$ by $(b_1\otimes b_2)^*= R(b_2^* \otimes b_1^*)$.
\end{remark}

\begin{remark} \label{rem:starclass}
For $\la\in\Pb^+$ define the involution $*$ on $B(\la)$ to be the
unique map that sends the highest weight vector $u_\la$ to the
lowest weight vector (the unique vector of weight $w_0(\la)$) and
satisfies \eqref{eq:defstar} for $i\in \J$. By \eqref{eq:sigalpha}
it follows that $\wt(b^*)=w_0\wt(b)$ for all $b\in B(\la)$.

Explicitly, the involution $*$ on the $U_q(\gehb)$-crystal
$B(\Lab_1)$ is given by
\begin{equation*}
\begin{aligned}[2]
& & i&\leftrightarrow \bar{i}  \\
& & \0&\leftrightarrow\0
\end{aligned}
\end{equation*}
except for
\begin{equation*}
\begin{aligned}[2]
\text{$\gehb=A_{n-1}$:}& &i&\leftrightarrow n+1-i \\
\text{$\gehb=D_n$, $n$ odd:} & & n&\leftrightarrow n \qquad
\bar{n}\leftrightarrow\bar{n}.
\end{aligned}
\end{equation*}
Here we use that the crystal of $B(\Lab_1)$ has underlying set~\cite{KN:1994}
\begin{align*}
  &\{1<2<\cdots<n\}
     && \text{for $A_{n-1}$}\\
  &\{1<2<\cdots<n<\0<\bar{n}<\cdots<\bar{2}<\bar{1} \}
     && \text{for $B_n$} \\
  &\{1<2<\cdots<n<\bar{n}<\cdots<\bar{2}<\bar{1} \}
     && \text{for $C_n$} \\
  &\{1<2<\cdots<\begin{matrix}{n} \\ {\bar{n}} \end{matrix} <\cdots<\bar{2}<\bar{1} \}
     && \text{for $D_n$.}
\end{align*}
\end{remark}

\subsection{Explicit formula for $*$}
We wish to determine the map $*$ of Conjecture
\ref{conj:crystalstar} explicitly for $B^s\in\CC$ and $B^{r,s}\in
\CCA$. The map $*:B^s\rightarrow B^s$ must stabilize classical
components by Remark \ref{rem:starcomp} and the
multiplicity-freeness of $B^s$ as a classical crystal. On each
classical component $B(s'\Lab_1)$ of $B^s$, $*$ is uniquely defined
by Remark \ref{rem:starclass}. Using the $U_q(\gehb)$-embedding
$B(s'\Lab_1)\rightarrow B(\Lab_1)^{\otimes s'}$ and Proposition
\ref{pp:star}, we have $(b_1b_2\dotsm b_{s'})^*=b_{s'}^* \dotsm
b_2^* b_1^*$. For $B^{r,s}\in \CCA$ and $b\in B^{r,s}$, $b^*$ is the
tableau obtained by replacing every entry $c$ of $b$ by $c^*$ and
then rotating by 180 degrees. The resulting tableau is sometimes
called the \textbf{antitableau} of $b$.

\begin{example}
For type $D_5^{(1)}$ we have
\Yboxdim{12pt}
$$
\young(113\be)^{\; *}=\young(\be\bc\ba\ba)
$$
For type $A_4^{(1)}$
\begin{equation*}
{\young(11,23)}^{\; *}=\young(34,55)\,.
\end{equation*}
\end{example}

\begin{prop} \label{pp:srstar} With $*$ defined as above,
Conjecture \ref{conj:crystalstar} holds for $B^s\in \CC$ and
$B^{r,s}\in \CCA$.
\end{prop}

\begin{remark} \label{rem:star}
{}From now on the notation $*$ will only be used in the following
way. Let $B=B_L\otimes\dotsm\otimes B_1$ be a tensor product of
factors $B_j=B^{s_j}$. Since $*$ may be regarded as an involution on
$B^s$, by Proposition \ref{pp:star} we may write
$B^*=B_1^*\otimes\dotsm\otimes B_L^*=B_1\otimes\dotsm\otimes B_L$
for the \textbf{reversed} tensor product. Then $*:B\rightarrow B^*$
is defined by $(b_L\otimes\dotsm\otimes b_1)^*\mapsto
b_1^*\otimes\dotsm\otimes b_L^*$.
\end{remark}

\begin{prop} \label{pp:R*} Let $R_j$ be the $R$-matrix acting at the $j^\nth$ and
$(j+1)^\fst$ tensor positions from the right. On an $L$-fold
tensor product of crystals of the form $B^s$,
\begin{equation} \label{eq:R*}
R_j \circ * = * \circ R_{L-j}
\end{equation}
for $1\le j\le L-1$.
\end{prop}
\begin{proof} One may reduce to the case $L=2$. Since $B_2\otimes B_1$ is connected,
$R$ is an isomorphism, and since \eqref{eq:defstar} holds, it suffices
to check \eqref{eq:R*} on $u(B_2\otimes B_1)$. But this holds by
weight considerations.
\end{proof}

\section{Virtual crystals}\label{s:virtual}
We review the virtual crystal construction \cite{OSS:2003a,OSS:2003b}. This allows
one to reduce the study of affine crystal graphs to those of simply-laced type.

\subsection{Embeddings of affine algebras} \label{ss:VX}
Any affine algebra $\geh$ of type $X$ can be embedded into a
simply-laced affine algebra $\geh_Y$ of type $Y$
\cite{HKOTT:2001}. For $\geh$ nonexceptional the embeddings are
listed below. The notation $A_{2n}^{(2)}$ and
$A_{2n}^{(2)\dagger}$ is used for two different vertex labelings
of the same Dynkin diagram, in which $\alpha_0$ is respectively
the extra short and extra long root.
\begin{equation} \label{eq:embed}
\begin{split}
  C_n^{(1)},A_{2n}^{(2)},A_{2n}^{(2)\dagger},D_{n+1}^{(2)} &\hookrightarrow
A_{2n-1}^{(1)} \\
  B_n^{(1)},A_{2n-1}^{(2)} &\hookrightarrow D_{n+1}^{(1)}.
\end{split}
\end{equation}

\subsection{Folding automorphism}
Let $\aut$ be the following automorphism of the Dynkin diagram of
$Y$. For $A_{2n-1}^{(1)}$, $\aut(i)=2n-i$ (mod $2n$). For type
$D_{n+1}^{(1)}$, $\aut$ exchanges the nodes $n$ and $n+1$ and
fixes all others.

Let $I^X$ and $I^Y$ be the vertex sets of the diagrams $X$ and $Y$
respectively, $I^Y/\aut$ the set of orbits of the action of $\aut$
on $I^Y$, and $\bij:I^X\rightarrow I^Y/\aut$ a bijection which
preserves edges and sends $0$ to $0$.

\begin{example}
If $Y=A_{2n-1}^{(1)}$, then $\bij(0)=\{0\}$, $\bij(i)=\{i,2n-i\}$
for $0<i<n$ and $\bij(n)=\{n\}$.

If $Y=D_{n+1}^{(1)}$, then $\bij(i)=i$ for $i<n$ and
$\bij(n)=\{n,n+1\}$.
\end{example}

\subsection{Embedding of weight lattices}
For $i\in I^X$ define $\mult_i$ as follows.
\begin{enumerate}
\item Let $Y=D_{n+1}^{(1)}$. 
\begin{enumerate}
\item Suppose the arrow points towards the component of $0$.
Then $\mult_i=1$ for all $i\in I^X$.
\item Suppose the arrow points away from the component of $0$.
Then $\mult_i$ is the order of $\aut$ for $i$ in the component of
$0$ and is $1$ otherwise.
\end{enumerate}
\item Let $Y=A_{2n-1}^{(1)}$. 
Then $\mult_i=1$ for $1\le i\le n-1$. For $i\in\{0,n\}$,
$\mult_i=2$ (which is the order of $\aut$) if the arrow incident
to $i$ points away from it and is $1$ otherwise.
\end{enumerate}

\begin{example}
For $X=B_n^{(1)}$ and $Y=D_{n+1}^{(1)}$ we have $\mult_i=2$ if
$0\le i \le n-1$ and $\mult_n=1$. For $X=A_{2n-1}^{(2)}$ and
$Y=D_{n+1}^{(1)}$ we have $\mult_i=1$ for all $i$.
\end{example}

The embedding $\emb:P^X\to P^Y$ of weight lattices is defined by
\begin{equation*}
\emb(\La^X_i) = \mult_i \sum_{j\in \bij(i)} \La^Y_j.
\end{equation*}
As a consequence we have
\begin{equation}\label{eq:Vdelta}
\begin{split}
  \emb(\alpha^X_i) &= \mult_i \sum_{j\in \bij(i)} \alpha^Y_j \\
  \emb(\delta^X) &= a_0^X \mult_0 \,\delta^Y.
\end{split}
\end{equation}

\subsection{Virtual crystals}
Fix an embedding $\geh_X\hookrightarrow \geh_Y$ in
\eqref{eq:embed}.

Let $\Vh$ be a $Y$-crystal. For $i\in I^X$ define the virtual
crystal operators $\eh_i,\fh_i$ on $\Vh$, as the composites of
$Y$-crystal operators $e_j,f_j$ given by
\begin{equation*}
\eh_i = \prod_{j\in \bij(i)} e_j^{\mult_i}\qquad\qquad \fh_i =
\prod_{j\in \bij(i)} f_j^{\mult_i} .
\end{equation*}

A \textbf{virtual crystal} (aligned in the sense of
\cite{OSS:2003a,OSS:2003b}) is an injection
$\emb:B\rightarrow \Vh$ from an $X$-crystal $B$ to a $Y$-crystal
$\Vh$ such that:
\begin{enumerate}
\item For all $b\in B$, $i\in I^X$, and $j\in \bij(i)\subset
I^Y$, $\vp_j(\emb(b))=\mult_i \vp_i(b)$ and
$\ve_j(\emb(b))=\mult_i \ve_i(b)$.
\item $\emb\circ e_i = \eh_i$ and $\emb\circ f_i = \fh_i$ for all $i\in
I^X$.
\end{enumerate}

A virtual crystal realizes the $X$-crystal $B$ as the subset of
the $Y$-crystal $\Vh$ given by its image under $\emb$, equipped
with the virtual Kashiwara operators $\eh_i$ and $\fh_i$.

A morphism $g$ of virtual crystals $\emb:B\rightarrow\Vh$ and
$\emb':B'\rightarrow \Vh'$ consists of a morphism
$g_X:B\rightarrow B'$ of $X$-crystals and a morphism
$g_Y:\Vh\rightarrow \Vh'$ of $Y$-crystals, such that the diagram
commutes:
\begin{equation*}
\begin{CD}
  B @>{\emb}>> \Vh \\
@V{g_X}VV @VV{g_Y}V \\
B' @>>{\emb'}> \Vh'
\end{CD}
\end{equation*}
An isomorphism $g$ of virtual crystals is a morphism $(g_X,g_Y)$
such that $g_X$ (resp. $g_Y$) is an isomorphism of $X$- (resp.
$Y$-) crystals.

\subsection{Tensor product of virtual crystals}
Let $\emb:B\rightarrow \Vh$ and $\emb':B'\rightarrow \Vh'$ be
virtual crystals. It is straightforward to verify that $\emb
\otimes \emb':B \otimes B' \rightarrow \Vh \otimes \Vh'$ is a
virtual crystal. Virtual crystals form a tensor category
\cite{OSS:2003a}.

\subsection{Virtual $B^s$}
We recall from \cite{OSS:2003b} the virtual crystal construction of
$B^s=B^{1,s}$ for $\geh$ of nonexceptional affine type. Let
$\Vh^s$ be given by
\begin{equation*}
  \Vh^s = \begin{cases}
B_Y^{s\vee} \otimes B_Y^s & \text{if $\geh_Y=A_{2n-1}^{(1)}$} \\
B_Y^s &\text{if $\geh_Y=D_{n+1}^{(1)}$ and $\geh=A_{2n-1}^{(2)}$} \\
B_Y^{2s} & \text{if $\geh_Y=D_{n+1}^{(1)}$ and $\geh=B_n^{(1)}$.}
\end{cases}
\end{equation*}

\begin{theorem} \label{th:V} \cite{OSS:2003b} There is a unique
virtual crystal $\emb:B^s\rightarrow \Vh^s$ such that
$\emb(u(B^s)) = u(\Vh^s)$.
\end{theorem}

\begin{example}
\Yboxdim{13pt} Let $X=B_3^{(1)}$ and $Y=D_4^{(1)}$. Then
$\Vh^s=B_Y^{2s}$. Let $b=\young(1\circ\bb)\in B_X^3$. Then
$\emb(b)=\young(113\bc\bb\bb)$ and $f_3(b)=\young(1\bc\bb)$.
Furthermore $$\fh_3(\emb(b))=f_3\circ
f_4(\emb(b))=\young(11\bc\bc\bb\bb).$$
\end{example}

\subsection{Virtual R-matrix}

\begin{prop} \label{pp:VR} \cite{OSS:2003b} Let
$\Rh:\Vh^t\otimes\Vh^s\rightarrow \Vh^s\otimes\Vh^t$ be the
composition of combinatorial $R$-matrices of type $Y$. Then the
diagram commutes:
\begin{equation*}
\begin{CD}
  B^t \otimes B^s @>{\emb}>> \Vh^t \otimes \Vh^s \\
@V{R}VV @VV{\Rh}V \\
B^s \otimes B^t @>>{\emb}> \Vh^s \otimes \Vh^t.
\end{CD}
\end{equation*}
That is, the pair $(R,\Rh)$ is an isomorphism of virtual crystals
$\emb:B^t\otimes B^s\rightarrow \Vh^t\otimes \Vh^s$ and
$\emb:B^s\otimes B^t\rightarrow \Vh^s \otimes \Vh^t$.
\end{prop}

\subsection{Virtual local coenergy}
\begin{prop} \label{pp:VH} \cite{OSS:2003b} Let $\emb:B\rightarrow\Vh$ and
$\emb':B'\rightarrow \Vh'$ be virtual crystals where $B,B'\in\CC$
both of type $X$. Then
\begin{equation*}
 H^X_{B,B'} = \frac{1}{\mult_0}\cdot H^Y_{\Vh,\Vh'} \circ (\emb\otimes
\emb').
\end{equation*}
\end{prop}

\subsection{Virtual graded crystal}
\begin{prop} \label{pp:VD} \cite{OSS:2003b} Let $B\in \CC$ be a crystal of type $X$ and
$\emb:B\rightarrow \Vh$ the corresponding virtual crystal. Then
\begin{equation*}
 D^X = \frac{1}{\mult_0}\cdot D^Y \circ \Psi.
\end{equation*}
\end{prop}

\subsection{Virtual $X$ formula}\label{ss:VX=V}
Let $B\in \CC$ be a crystal of type $X$. Let $\emb:B\rightarrow\Vh$ be the corresponding
virtual crystal. For $\la\in\Pb^+$ let $\Pv(B,\la)$ be image under $\emb$ of
the set $P(B,\la)$. Define the virtual X formula by
\begin{equation*}
  VX_{B,\la}(q) = \sum_{b\in \Pv(B,\la)} q^{D(\emb(b))/\mult_0}.
\end{equation*}

\begin{theorem} \label{th:X=VX} ($X=VX$) \cite{OSS:2003b} For $\geh$ of nonexceptional
affine type $X$ and $B\in\CC$ a crystal of type $X$, one has $X_{B,\la}(q)=VX_{B,\la}(q)$.
\end{theorem}

\subsection{Virtual crystals and $*$-duality}
We believe that the following is true for any virtual crystal,
namely, that up to $R$-matrices, $\emb$ takes the $*$ involution of
type $X$ to the $*$ involution of type $Y$.

\begin{prop} \label{pp:*emb} Let $B^s\in\CC$ and let
$\emb:B^s\rightarrow\Vh^s$ be a virtual crystal. Then the following
diagram commutes, where $\iota$ is either a composition of
$R$-matrices or the identity:
\begin{equation*}
\begin{CD}
B^s @>{\emb}>> \Vh^s \\
@V{*}VV @VV{\iota\circ *}V \\
B^s @>>{\emb}> \Vh^s
\end{CD}
\end{equation*}
\end{prop}
\begin{proof} Note that for the non-simply-laced types $X$ the Dynkin
involution $\tau_X$ is the identity. The virtual raising and
lowering operators are invariant under $\tau_Y$. It is therefore
sufficient to check the above commutation on $v\in P(B^s)$, where
$P(B^s)$ is the set of classical highest weight vectors in $B^s$. But
$B^s$ and $\Vh^s$ are multiplicity-free as a classical crystals and
$*$ stabilizes classical components and modifies the weight of a
crystal element by applying $w_0$. The following are equivalent:
\begin{enumerate}
\item $\emb(v)^*\in \Vh^s$ is a classical lowest weight vector and
$\wt(\emb(v)^*)=w_0^Y(\emb(\la))$.
\item $\emb(v)\in P(\Vh^s)$ and $\wt(\emb(v))=\emb(\la)$.
\item $v\in P(B^s)$ and $\wt(v)=\la$.
\item $v^*\in B^s$ is a classical lowest weight vector and $\wt(v^*)=w_0^X\la$.
\item $\emb(v^*)$ is a classical lowest weight vector and
$\wt(\emb(v^*))=\emb(w_0^X\la)$.
\end{enumerate}
One may verify that $w_0^Y(\emb(\la))=\emb(w_0^X(\la))$ using
linearity, to reduce to the case $\la=\Lab_i^X$ for $i\in\J$. It
follows that $\emb(v)^*$ and $\emb(v^*)$ are classical lowest weight
vectors in $\Vh^s$ of the same weight. But then they must be equal.
\end{proof}

\section{Right splitting}\label{s:right}
Let $\geh$ be of nonexceptional affine type. We define a family of
$U_q(\gehb)$-crystal embeddings which is well-behaved with respect
to intrinsic coenergy. They are denoted $\rs:=\rs_{r;a,b}$ which stands for
``right-split", because when $b=0$, the map splits off the
rightmost column of an element in $B^{r,s}$.

\begin{conjecture} \label{conj:split} Let $a-2\ge b\ge0$. Suppose
$\CC'$ is a set of KR crystals whose modules have been constructed,
which contains $B^{r,s}$ for a particular $r\in\J$ and all
$s\in\Z_{>0}$. Then there is an injective $U_q(\gehb)$-crystal
morphism
\begin{equation*}
 \rs_{r;a,b}: B^{r,a} \otimes B^{r,b} \rightarrow
B^{r,a-1}\otimes B^{r,b+1}
\end{equation*}
such that for any crystal $B$ which is the tensor product of
crystals in $\CC'$, the map
\begin{equation} \label{eq:1rsmap}
  1_B \otimes \rs_{r;a,b}: B \otimes B^{r,a} \otimes B^{r,b}
\rightarrow B \otimes B^{r,a-1} \otimes B^{r,b+1}
\end{equation}
is an injective $U_q(\gehb)$-crystal morphism which preserves
intrinsic coenergy.
\end{conjecture}

\begin{theorem} \label{th:Asplit} Conjecture \ref{conj:split}
holds for $\geh=A_n^{(1)}$ for all $r\in\J$ and the set $\CC'$ of
all KR crystals.
\end{theorem}
\begin{proof} This follows from \cite{Sh:2001a,Sh:2001b,Sh:2002}.
\end{proof}

\begin{theorem} \label{th:splitrow} Conjecture \ref{conj:split}
holds for any nonexceptional affine algebra $\geh$ for $r=1$ and
$\CC'$ the set of KR crystals of the form $B^{1,s}$.
\end{theorem}
The proof of Theorem \ref{th:splitrow} occupies the remainder of
this section.

\subsection{Explicit definition of splitting}
This paper only requires the case $b=0$ for the map $\rs$. Except
for type $A_n^{(1)}$ only the case $r=1$ is needed. For $s\ge2$
define the map $\rs:=\rs_{1;s,0}$ as follows. For types $A_n^{(1)}$,
$D_n^{(1)}$, $B_n^{(1)}$, and $A_{2n-1}^{(2)}$, define
$\rs:B^s\rightarrow B^{s-1} \otimes B^1$ by $\rs(wx)=w\otimes x$ for
$x\in B^1$ and $w\in B^{s-1}$ such that $wx\in B^s$. For the other
types, in addition to the above rules we have $\rs(x)=\vn\otimes x$
for $x\in B(\Lab_1)\subseteq B^s$, and $\rs(\vn)=\bar{1}\otimes 1$.
For $B^{r,s}\in \CCA$ and $b\in B^{r,s}$, let $\rs(b)=b_2\otimes
b_1$, where $b_1$ is the rightmost column of the rectangular tableau
$b$ and $b_2$ is the rest of $b$.

\begin{remark} \label{rem:rs} Suppose $s\ge 2$. Here $r=1$
for $\CC$. For $B^{r,s}\in \CC$ (or $\CCA$) we write $\rs$ for the
map $1_B \otimes \rs$ on $B \otimes B^{r,s}$ and write $\rs(B\otimes
B^{r,s}):=B\otimes B^{r,s-1}\otimes B^{r,1}$.
\end{remark}

\subsection{Simply-laced $\geh$}
$\geh=A_n^{(1)}$ is covered by Theorem \ref{th:Asplit}. The other
simply-laced nonexceptional family is $\geh=D_n^{(1)}$.

It is straightforward to check directly using the explicit
description of $B^s$ in \cite{OSS:2003b} that $\rs$ is an injective
$U_q(\gehb)$-crystal morphism. Let $B$ be the tensor product of
crystals in $\CC'$. To check that $1_B\otimes \rs$ preserves
intrinsic coenergy, by \eqref{eq:DNY} it suffices to check this
property for $B$ the trivial crystal and for $B=B^t$. Since $1_B
\otimes \rs$ is a $U_q(\gehb)$-crystal morphism, it is sufficient to
prove that intrinsic coenergy is preserved for classical highest
weight vectors. Suppose $B$ is trivial. By \eqref{eq:ClDecomp} $B^s$
has a single classical highest weight vector, namely, $u(B^s)=1^s$.
By Example \ref{ex:DKR} $D_{B^s}=0$. On the other hand
$\rs(1^s)=1^{s-1} \otimes 1 = u(B^{s-1} \otimes B^1)$ so its
intrinsic energy is also zero. For $B=B^t$ we require the following
Lemma, which is easily verified directly.

\begin{lemma} \label{lem:D2P} For $\geh=D_n^{(1)}$ and $s,t\ge 1$,
$P(B^t\otimes B^s)$ consists of the elements
\begin{equation*}
  v_{p,q}^{t,s} = 1^{t-p-q}\,2^p\,\bar{1}^q \otimes 1^s
\end{equation*}
where $p+q\le\min(s,t)$. In particular $B^t \otimes B^s$ is
multiplicity-free as a $U_q(D_n)$-crystal.
\end{lemma}

Recall that $D_{B^t}$ and $D_{B^s}$ are identically zero by Example
\ref{ex:DKR}. By \eqref{eq:DNY} and explicit calculation,
\begin{equation} \label{eq:DD2}
  D_{B^t\otimes B^s}(v_{p,q}^{t,s}) = H_{B^t,B^s}(v_{p,q}^{t,s})=p+2q.
\end{equation}
Since $R$ is a $U_q(D_n^{(1)})$-crystal isomorphism, Lemma
\ref{lem:D2P} implies that
\begin{equation*}
  R_{B^t,B^s}(v_{p,q}^{t,s}) = v_{p,q}^{s,t}.
\end{equation*}
We compute $H_{B^t,B^{s-1} \otimes B^1}$ on
$\rs(v_{p,q}^{t,s})=(1^{t-p-q}2^p\bar{1}^q\otimes 1^{s-1} \otimes
1)$ using Proposition \ref{pp:Hbraided}. We have
\begin{equation*}
R_{B^t,B^{s-1}}(1^{t-p-q}2^p\bar{1}^q \otimes 1^{s-1})=
\begin{cases}
1^{s-1-p-q}2^p\bar{1}^q \otimes 1^t
 & \text{if $p+q<s$}\\
\overline{1}^{q-1} \otimes 1^{t-1} \bar{1}
 & \text{if $p+q=s$, $q=s$}\\
2^{p-1} \overline{1}^q \otimes 1^{t-1}2
 & \text{if $p+q=s$, $q<s$.}
\end{cases}
\end{equation*}
By \eqref{eq:DD2} we have
\begin{equation*}
H_{B^t,B^{s-1}}(1^{t-p-q}2^p\bar{1}^q \otimes 1^{s-1})=
\begin{cases}
p+2q & \text{if $p+q<s$}\\
p+2q-2 & \text{if $p+q=s$, $q=s$}\\
p+2q-1 & \text{if $p+q=s$, $q<s$.}
\end{cases}
\end{equation*}
By \eqref{eq:DD2} we have $H(1^t\otimes 1)=0$,
$H(1^{t-1}\bar{1}\otimes 1)=2$, and $H(1^{t-1}2\otimes 1)=1$. It
follows in any case that $1_{B^t}\otimes \rs$ preserves intrinsic
coenergy.

\subsection{Non-simply-laced $\geh$}
Suppose $\geh$ is not simply-laced. Let $\geh\hookrightarrow \geh_Y$
be as in \eqref{eq:embed}. It is not hard to show that $\rs_X$ is a
$U_q(\gehb)$-crystal injection. To show that the map
\eqref{eq:1rsmap} preserves intrinsic coenergy (and thereby complete
the proof of Theorem \ref{th:splitrow}), by Proposition \ref{pp:VD}
the following result suffices.

\begin{prop} \label{pp:vrs} There is an injective $U_q(\geh_Y)$-crystal
map $\vrs:\Vh^s\rightarrow \Vh^{s-1} \otimes \Vh^1$ such that:
\begin{enumerate}
\item \label{it:vrsdiag} The following diagram commutes:
\begin{equation} \label{eq:vrsdiag}
\begin{CD}
B^s_X @>{\Psi}>> \Vh^s\\
@V{\rs_X}VV @VV{\vrs}V \\
B_X^{s-1} \otimes B_X^1  @>>{\Psi\otimes\Psi}> \Vh^{s-1} \otimes
\Vh^1
\end{CD}
\end{equation}
\item \label{it:vD} For any $B\in\CC$, let $\emb:B\rightarrow\Vh$ be its virtual
crystal embedding. Then $1_{\Vh} \otimes \vrs$ preserves intrinsic
coenergy.
\item \label{it:vres} If $v\in \Vh^s$ and
$\vrs(v)\in\Image(\Psi\otimes\Psi)$ then $v\in\Image(\Psi)$.
\end{enumerate}
\end{prop}
\begin{proof} Suppose first that $Y=A_{2n-1}^{(1)}$. Then $\Vh^s =
B_Y^{s\vee} \otimes B_Y^s$. Define the map $\vrs:\Vh^s\rightarrow
\Vh^{s-1}\otimes \Vh^1$ by the composition
\begin{equation} \label{eq:vrsAdef}
\begin{split}
 &B_Y^{s\vee} \otimes B_Y^s \stackrel{1\otimes \rs_Y}{\longrightarrow}
B_Y^{s\vee} \otimes B_Y^{s-1} \otimes B_Y^1\stackrel{R}{\longrightarrow}
B_Y^{s-1} \otimes B_Y^1 \otimes B_Y^{s\vee}\\
\stackrel{1\otimes 1\otimes \rs_Y^\vee}{\longrightarrow}&
B_Y^{s-1} \otimes B_Y^1 \otimes B_Y^{s-1\vee} \otimes B_Y^{1\vee}
\stackrel{R}{\longrightarrow}
B_Y^{s-1\vee}\otimes B_Y^{s-1} \otimes B_Y^{1\vee} \otimes B_Y^1.
\end{split}
\end{equation}
Here $\rs^\vee_Y(wx)=w\otimes x$ where $wx\in B_Y^{s\vee}$, $w\in
B_Y^{s-1\vee}$ and $x\in B_Y^{1\vee}$. Note that \eqref{eq:vrsAdef}
is a composition of combinatorial $R$-matrices and $\rs$ maps for
type $A$. Point \ref{it:vD} holds by Theorem \ref{th:Asplit}.

One need only verify \ref{it:vrsdiag} on classical highest weight
vectors, by the definition of $\Psi$ and the fact that $\rs_Y$ and
$\rs^\vee_Y$ (resp. $\rs_X$) are morphisms of $U_q(\overline{Y})$-
(resp. $U_q(\overline{X})$-) crystals.

Let $N=2n$. The classical highest weight vectors in $B_X^s$ have the
form $1^{s-p}$ for $0\le p\le s$; if $\geh$ is $C_n^{(1)}$ or
$A_{2n}^{(2)\dagger}$ then $p$ must also be even. For $p<s$ the
element $1^{s-p}\in B((s-p)\Lab_1)\subset B_X^s$ is sent to the
following elements under the maps in \eqref{eq:vrsdiag}:
\begin{equation*}
\begin{matrix}
1^{s-p} & 1^{s-1-p} \otimes 1  \\ N^{\vee(s-p)} 1^{\vee p} \otimes
1^s & N^{\vee s-p-1} 1^{\vee p} \otimes 1^{s-1} \otimes N^\vee
\otimes 1
\end{matrix}
\end{equation*}
where the intermediate results under the maps in
\eqref{eq:vrsAdef} are given by $N^{\vee(s-p)} 1^{\vee p} \otimes
1^{s-1} \otimes 1$, $1^{s-p-1}N^p \otimes 1 \otimes N^{\vee s}$,
$1^{s-p-1}N^p \otimes 1 \otimes N^{\vee(s-1)} \otimes N^\vee$, and
$N^{\vee s-p-1} 1^{\vee p} \otimes 1^{s-1} \otimes N^\vee \otimes
1$.

Under the maps in \eqref{eq:vrsdiag}, the element $\es$ is sent to
\begin{equation*}
\begin{matrix}
\es & \overline{1} \otimes 1 \\
N^s \otimes N^{\vee s} & 1^{\vee s-1} \otimes 1^{s-2}N \otimes
N^\vee \otimes 1
\end{matrix}
\end{equation*}
with intermediate values in \eqref{eq:vrsAdef} given by $N^s
\otimes N^{\vee s-1} \otimes N^\vee$, $1^{\vee s-1} \otimes 1^\vee
\otimes 1^s$, $1^{\vee s-1} \otimes 1^\vee \otimes 1^{s-1} \otimes
1$, and $1^{\vee s-1} \otimes 1^{s-2}N \otimes N^\vee \otimes 1$.

Since these are all the possible classical highest weight vectors,
point \ref{it:vrsdiag} follows.

For point \ref{it:vres}, let $v\in \Vh^s$ and
$\vrs(v)\in\Image(\emb\otimes\emb)$. Without loss of generality we
may assume that $v\in P(\Vh^s)$ since $\vrs$ is a
$U_q(\overline{Y})$-morphism. Now $v$ must have the form
$v_{s,p}:=N^{\vee(s-p)} 1^{\vee p} \otimes 1^s$ for $0\le p\le s$.
By computations similar to those above, $\vrs(v)=v_{s-1,p} \otimes
\emb(1)$ if $p<s$ and $\vrs(v)=1^{\vee(s-1)}\otimes 1^{s-2}N\otimes
\emb(1)$ if $p=s$. But $\vrs(v)\in\Image(\emb\otimes\emb)$ means
that $v_{s-1,p}\in \Image(\emb)$ if $p<s$ and $1^{\vee(s-1)}\otimes
1^{s-2}N\in\Image(\emb)$ if $p=s$. The parity condition for this to
occur implies the parity condition that guarantees that
$v_{s,p}\in\Image(\emb)$.

Suppose next that $Y=D_{n+1}^{(1)}$ and $X=A_{2n-1}^{(2)}$. Then
$\Vh^s=B_Y^s$. Define $\vrs=\rs_Y:B_Y^s\rightarrow B_Y^{s-1}\otimes
B_Y^1$. Point \ref{it:vD} follows by the simply-laced $D_n^{(1)}$
case. Point \ref{it:vres} is trivial. For point \ref{it:vrsdiag} it
is enough to consider elements of $P(B^s)=\{1^s\}$. Under the maps
in \eqref{eq:vrsdiag}, $1^s$ goes to
\begin{equation*}
\begin{matrix}
1^s & 1^{s-1} \otimes 1 \\
1^s & 1^{s-1} \otimes 1
\end{matrix}
\end{equation*}
and \eqref{eq:vrsdiag} commutes.

Suppose that $Y=D_{n+1}^{(1)}$ and $X=B_n^{(1)}$. Then
$\Vh^s=B_Y^{2s}$. Define $\vrs:B_Y^{2s} \rightarrow B_Y^{2s-2}
\otimes B_Y^2$ by $wv\mapsto w\otimes v$ where $wv\in B_Y^{2s}$,
$w\in B_Y^{2s-2}$, and $v\in B_Y^2$. This map is clearly injective
and $U_q(\overline{Y})$-equivariant. Point \ref{it:vres} is obvious.
For point \ref{it:vrsdiag} it is enough to consider the unique
element $1^s\in P(B_X^s)$. Under \eqref{eq:vrsdiag} $1^s$ goes to
\begin{equation*}
\begin{matrix}
1^s & 1^{s-1} \otimes 1 \\
1^{2s} & 1^{2s-2} \otimes 1^2
\end{matrix}
\end{equation*}
so that \eqref{eq:vrsdiag} commutes. For point \ref{it:vD} define
$\vrs':B_Y^{2s}\rightarrow B_Y^{2s-2} \otimes B_Y^1 \otimes B_Y^1$
by the composite map
\begin{equation*}
\begin{split}
B_Y^{2s} &\stackrel{\rs_Y}{\longrightarrow} B_Y^{2s-1}\otimes B_Y^1
\stackrel{R}{\longrightarrow} B_Y^1\otimes B_Y^{2s-1}\\
&\stackrel{1\otimes \rs_Y}{\longrightarrow} B_Y^1\otimes B_Y^{2s-2}\otimes B_Y^1
\stackrel{R}{\longrightarrow} B_Y^{2s-2} \otimes B_Y^1 \otimes B_Y^1.
\end{split}
\end{equation*}
Since $\vrs'$ is the composition of $\rs_Y$ maps and $R$-matrices,
it preserves intrinsic coenergy by the simply-laced case. It
suffices to show that
\begin{equation*}
\xymatrix{
   B_Y^{2s} \ar[dd]^{\vrs} \ar[dr]^{\vrs'} & \\
  {} & {B_Y^{2s-2} \otimes B_Y^1 \otimes B_Y^1} \\
   {B_Y^{2s-2} \otimes B_Y^2} \ar[ur]_{1\otimes \rs_Y} &
}
\end{equation*}
commutes since $\vrs'$ and $1\otimes \rs_Y$ both preserve intrinsic
coenergy. It suffices to check this for the lone classical highest
weight vector $1^{2s}\in P(B_Y^{2s})$. Clearly
$\vrs'(1^{2s})=1^{2s-2}\otimes 1\otimes1$, while
$\vrs(1^{2s})=1^{2s-2}\otimes 1^2$ and this is sent by $1\otimes
\rs_Y$ to $1^{2s-2}\otimes 1\otimes 1$, as desired.
\end{proof}

\section{Left splitting and duality}\label{s:left}
We define dual analogues of the intrinsic coenergy $D$ and right
splitting.

\subsection{Tail coenergy} For $B^s\in \CC$ define
$\Dt_{B^s}=D_{B^s}$. For $B^{r,s}\in \CCA$, define
$\Dt_{B^{r,s}}=D_{B^{r,s}}=0$. If $B_1,B_2,\dotsc,B_L\in \CC$ (or
$\CCA$) and $B=B_L\otimes\dotsm\otimes B_1$ are such that
$\Dt_{B_j}:B_j\rightarrow\Z_{\ge0}$ are given, then define
\begin{equation} \label{eq:tailDNY}
  \Dt_B = \sum_{1\le i<j\le L} H_{j-1} R_{j-2} \dotsm R_{i+1} R_i
   + \sum_{j=1}^L \Dt_{B_j} R_{L-1} R_{L-2} \dotsm R_j
\end{equation}
with $\Dt_{B_j}$ acting on the leftmost tensor position. This is a
different associative tensor product on graded crystals than the
one given in subsection \ref{ss:D}.

Recall the notation $B^*$ of Remark \ref{rem:star}.

\begin{prop} \label{pp:tail} Let $B\in \CC$ (or $B\in \CCA$) and $b\in B$.
Then $\Dt_B(b)=D_{B^*}(b^*)$.
\end{prop}
\begin{proof} For $B$ a single KR crystal, the result follows from
the fact that the involution $*$ on $B$ stabilizes classical
components. By Proposition \ref{pp:R*} and comparing
\eqref{eq:tailDNY} with \eqref{eq:DNY} it suffices to show that
\begin{equation}\label{eq:H*}
H_{B_1,B_2}(b^*)=H_{B_2,B_1}(b)
\end{equation}
for $B_1,B_2$ KR crystals. Since $B_2\otimes B_1$ is connected, the
proof may proceed by induction on the number of steps (either of the
form $e_i$ or $f_i$) in $B_2\otimes B_1$ from $u(B_2\otimes B_1)$ to
$b$. Suppose first that $b=u(B_2\otimes B_1)$. By the definition of
$u(B)$ in subsection \ref{ss:fin}, $B_2 \otimes B_1$ (and therefore
$B_1\otimes B_2)$ contain a unique classical component isomorphic to
$B(\la)$ where $\la=\wt(b)$. And $B(\la)$ contains a unique vector
of the extremal weight $w_0\la$. Since $\wt(b^*)=w_0\wt(b)$ it
follows that $b^*$ and $u(B_1\otimes B_2)$ are in the same classical
component, so that $H_{B_1,B_2}(b^*)=H_{B_1,B_2}(b)=0$ by the
definition of $H$.

Now suppose $b=f_i(c)$ where $c$ is closer to $u(B_2\otimes B_1)$
than $b$ is. If $i\not=0$ then we are done since both sides of
\eqref{eq:H*} do not change under passing from $c$ to $b$, by the
definition of $H$ and \eqref{eq:defstar}. So assume $i=0$. By
\eqref{eq:defstar} $b^*=e_0(c^*)$. But then one may conclude the
validity of \eqref{eq:H*} for $b$ from that of $c$ using rules for
the Kashiwara operators on the tensor product and \eqref{eq:H}.
\end{proof}

Define $\Xt$ just like the one-dimensional sum $X$ but use $\Dt_B$
instead of $D_B$. Proposition \ref{pp:tail} has this corollary.

\begin{corollary} $\Xt(B,\la)=X(B,\la)$.
\end{corollary}

\subsection{Left splitting} Whenever the right splitting map
$\rs:B^{r,s}\rightarrow B^{r,s-1} \otimes B^{r,1}$ is defined, we
may define the left-splitting map $\ls:B^{r,s}\rightarrow
B^{r,1}\otimes B^{r,s-1}$ by the commutation of the diagram
\begin{equation} \label{eq:splitstar}
\begin{CD}
  B^{r,s} @>{\ls}>> B^{r,1} \otimes B^{r,s-1} \\
  @V{*}VV @VV{*}V \\
  B^{r,s} @>>{\rs}> B^{r,s-1} \otimes B^{r,1}.
\end{CD}
\end{equation}
In particular, it is defined for $B^s\in \CC$ and $B^{r,s}\in\CCA$.

\begin{corollary} \label{cor:lsplitdef} Here $r=1$ for the category $\CC$.
$\ls$ is a $U_q(\gehb)$-crystal embedding such that, for any
$B^{r,s}\in \CC$ (or $\CCA$) and for any $B\in \CC$ (or $\CCA$), the map
\begin{equation*}
\ls \otimes 1_B: B^{r,s} \otimes B \rightarrow B^{r,1} \otimes
B^{r,s-1}\otimes B
\end{equation*}
is injective and preserves $\Dt$.
\end{corollary}
\begin{proof} $\ls$ is a $U_q(\gehb)$-crystal embedding since
$\rs$ is, by Theorem \ref{th:splitrow}, the definition of $*$ and
\eqref{eq:splitstar}. For the preservation of $\Dt$, let $b_1\otimes
b_2\in B^{r,s}\otimes B$. We have
\begin{equation*}
\Dt(\ls(b_1) \otimes b_2) =
D(b_2^* \otimes \ls(b_1)^*) =
D(b_2^* \otimes \rs(b_1^*)) =
D(b_2^* \otimes b_1^*) =
\Dt(b_1 \otimes b_2)
\end{equation*}
by Proposition \ref{pp:tail} and \eqref{eq:splitstar}.
\end{proof}

\begin{remark} \label{rem:splitnote} Suppose $s\ge 2$. Here $r=1$
for $\CC$ as usual. For $B^{r,s}\in \CC$ (or $\CCA$) we write $\ls$
for the map $\ls \otimes 1_B$ on $B^{r,s}\otimes B$. Also we write
$\ls(B^{r,s}\otimes B) := B^{r,1} \otimes B^{r,s-1}\otimes B$.
\end{remark}

\subsection{Explicit left-splitting}

\begin{lemma} \label{lem:lsplitrow} For $B^s\in \CC$ the map
$\ls:B^s\rightarrow B^1 \otimes B^{s-1}$ is given explicitly by
$\ls(xw)=x\otimes w$ for $x\in B^1$ and $w\in B^{s-1}$ such that
$xw\in B^s$, $\ls(x)=x\otimes \vn$ for $x\in B(\Lab_1)\subseteq
B^s$, $\ls(\vn)=\bar{1}\otimes 1$. For $B^{r,s}\in \CCA$ and $b\in
B^{r,s}$, $\ls(b)=b_2 \otimes b_1$ where $b_2$ is the leftmost
column in the $r\times s$ semistandard tableau $b$ and $b_1$ is the
rest of $b$.
\end{lemma}

\subsection{Box-splitting} \label{ss:bs}
Let $B^{r,1}\in \CCA$ with $r\ge2$. There is a $U_q(\gehb)$-crystal
embedding $\llb:B^{r,1}\rightarrow B^{1,1} \otimes B^{r-1,1}$ given
by $b\mapsto b_2 \otimes b_1$ where $b_2$ is the bottommost entry in
the column tableau $b$ of height $r$, and $b_1$ is the remainder of
$b$. There is a $U_q(\gehb)$-crystal embedding
$\rb:B^{r,1}\rightarrow B^{r-1,1}\otimes B^{1,1}$ given by $b\mapsto
b_2 \otimes b_1$ where $b_1$ is the topmost entry in the column $b$
and $b_2$ is the rest of $b$.

The map $\llb$ is only used to define the path-RC bijection for
$B\in\CCA$ in section~\ref{sec:bij}.

In general the morphism $\rb$ does not preserve intrinsic coenergy,
but another grading called intrinsic energy. It was proved in
\cite{KSS:2002} that the path-RC bijection preserves the grading for
$\CCA$ using a different method, namely, the rank-level duality for
type $A^{(1)}$.

\subsection{Projections and commutations}
Define the (``left-hat") map $\lh:B_2\otimes B_1\rightarrow B_1$ by
$b_2\otimes b_1\mapsto b_1$. It just removes the left tensor factor.
Define the ``right-hat" map $\rh:B_2\otimes B_1\rightarrow B_2$ by
$b_2\otimes b_1\mapsto b_2$.

It is immediate that the following diagram commutes:
\begin{equation} \label{eq:hat*}
\begin{CD}
  B_2 \otimes B_1 @>{\lh}>> B_1 \\
  @V{*}VV @VV{*}V \\
  B_1 \otimes B_2 @>>{\rh}> B_1
\end{CD}
\end{equation}

Let $P(B)$ be the set of classical highest weight vectors in $B$,
or equivalently, the set of classical components of $B$.

\begin{lemma} \label{lem:hwvhat}
The maps $\lh:B_2\otimes B_1\rightarrow B_1$ and $\rh:B_2\otimes
B_1\rightarrow B_2$ induce maps $\lh:P(B_2\otimes B_1)\rightarrow
P(B_1)$ and $\rh:P(B_2\otimes B_1) \rightarrow P(B_2)$.
\end{lemma}
\begin{proof} If $b_2\otimes b_1$ is a classical highest weight
vector of $B_2\otimes B_1$ then by the definitions, $b_1$ is a
classical highest weight vector of $B_1$. Thus $\lh$ is well-defined
on components.

For $\rh$ we work with classical components. By \eqref{eq:defstar}
the map $*$ takes classical components to classical components.
But then $\rh$ is well-defined on components since $\lh$ is, by
\eqref{eq:hat*}.
\end{proof}

\begin{example}
\Yboxdim{13pt} Let $b=\young(3) \otimes \young(2\bb) \otimes
\young(12) \otimes \young(1)\in P(B^1\otimes B^2\otimes B^2 \otimes
B^1)$ of type $D_4^{(1)}$. Then $\lh(b)=\young(2\bb) \otimes
\young(12) \otimes \young(1)$ and $\rh(b)=\young(3) \otimes
\young(2\bb) \otimes \young(12)$. The induced map on highest weight
vectors yields
$\rh(b)=\young(3)\otimes\young(22)\otimes\young(11)\;$
\end{example}

One has the commutation of induced maps on classical highest weight
vectors:
\begin{equation*}
\begin{CD}
  P(B_2 \otimes B_1) @>{\lh}>> P(B_1) \\
  @V{*}VV @VV{*}V \\
  P(B_1 \otimes B_2) @>>{\rh}> P(B_1)
\end{CD}
\end{equation*}

\begin{remark} \label{rem:hatnote}
{}From now on, unless explicitly indicated otherwise, we only
consider the map $\lh$ (resp. $\rh$) on tensor products whose left
(resp. right) factor is $B^1$. In these cases, we use the notation
$\lh(B^1\otimes B)=B$ and $\rh(B\otimes B^1)=B$.
\end{remark}
For $\la\in\Pb^+$ let
\begin{equation*}
  \lm = \{ \mu\in \Pb^+ \mid \text{$B(\la)$ occurs in $B^1 \otimes B(\mu)$ } \}
\end{equation*}
where $B^1$ is regarded as a $U_q(\gehb)$-crystal by restriction.

By Lemma \ref{lem:hwvhat} there are well-defined bijections
\begin{align}
\label{eq:lhbij}
\lh:P(B,\la)&\rightarrow \bigcup_{\mu\in \lm} P(\lh(B),\mu) \\
\label{eq:rhbij}
\rh:P(B,\la)&\rightarrow \bigcup_{\mu\in \lm} P(\rh(B),\mu)
\end{align}
except in the case $\geh=D_{n+1}^{(2)}$. Note that $B^1$ has at
most one vector of each weight except when $\geh=D_{n+1}^{(2)}$,
which has two vectors $0$ and $\vn$ of weight $0$. If $\mu=\la$,
then there can be elements $b\in P(\lh(B),\la)$ such that both
$0\otimes b$ and $\vn\otimes b$ are in $P(B,\la)$. If so, then the
right hand side of \eqref{eq:lhbij} must be modified to include
two copies of $b$, one coming from $\vn\otimes b$ and the other
from $0\otimes b$. There is no analogous problem for $\rh$ since
$0\not\in P(B^1)$.

\begin{prop} \label{pp:cryscomm}\mbox{} Let $r=r'=1$ for $\CC$.
\begin{enumerate}
\item \label{pt:lhrh} $[\lh,\rh]=0$ on $B^1 \otimes B \otimes B^1$.
\item \label{pt:lhrs} $[\lh,\rs]=0$ on $B^1 \otimes B \otimes
B^{r,s}$ for $s\ge2$.
\item \label{pt:rhls} $[\rh,\ls]=0$ on $B^{r,s} \otimes B \otimes
B^1$ for $s\ge2$.
\item \label{pt:lsrs} $[\ls,\rs]=0$ on $B^{r,s} \otimes B \otimes
B^{r',s'}$ for $s,s'\ge2$.
\item \label{pt:hat*} $*\circ \lh = \rh \circ *$ on $B^1 \otimes B$.
\end{enumerate}
Moreover, these commutations also hold for the induced maps on sets
of classical highest weight vectors.
\end{prop}
\begin{proof} The operators on the entire crystals commute more or
less by definition. We now prove that these identities hold for the
induced maps between sets of classical highest weight vectors.

The proof is again trivial except for cases involving $\rh$. Point
\ref{pt:lhrh} follows from Lemma \ref{lem:hwvhat}. Point
\ref{pt:rhls} follows from Lemma \ref{lem:hwvhat} and the
$U_q(\gehb)$-equivariance of $\ls$ given in Corollary
\ref{cor:lsplitdef}. Finally, point \ref{pt:hat*} follows from Lemma
\ref{lem:hwvhat} and the fact \eqref{eq:defstar} that the map $*$
respects classical raising and lowering operators.
\end{proof}

\subsection{Right hat and classical highest weight vectors}
We need to know precisely how the highest weights change when
passing from an element of $P(B^1\otimes B \otimes B^1)$ to $P(B)$
via either $\rh \circ \lh$ or $\lh \circ \rh$. In this section we
assume type $D_n$. The answer is given by van Leeuwen
\cite{vL:1998}. We translate his answer into the language of
partitions. 

Let $P$ be the set of dominant weights that can occur in a tensor
product of crystals of the form $B(\Lab_1)$. A dominant weight
$\sum_{i=1}^n a_i \Lab_i$ is in $P$ if and only if $a_{n-1}$ and
$a_n$ have the same parity. We put a graph structure on $P$ by
declaring that weights $\la$ and $\mu$ are adjacent if there is an
element $x\in B(\Lab_1)$ such that $\la-\mu=\wt(x)$.

We realize $P$ as a subset of $\Z^n$ by letting
$\Lab_i=(1^i,0^{n-i})$ for $1\le i\le n-2$,
$\Lab_{n-1}=\frac{1}{2}(1^n)$ and $\Lab_n=\frac{1}{2}(1^{n-1},-1)$.
As such $P$ is given by the tuples
$\la=(\la_1,\la_2,\dotsc,\la_n)\in\Z^n$ with $\la_1\ge \la_2
\ge\dotsm \ge \la_{n-1} \ge | \la_n |$.

We modify this notation slightly in order to use partitions. Let $Y$
be the lattice of partitions
$\la=(\la_1\ge\la_2\ge\dotsm\ge\la_n)\in\Z^n_{\ge0}$ with at most
$n$ parts. A graph structure on $Y$ is given by declaring that two
partitions are connected with an edge if their partition diagrams
differ by one cell. Define the graph $G$ by glueing two copies $Y_+$
and $Y_-$ of $Y$ together such that, if $\la\in Y$ is such that
$\la_n=0$, then $\la \in Y_+$ and $\la\in Y_-$ are identified.

Then $P \cong G$ where the weight $(\mu_1,\dotsc,\mu_n)$ is
identified with the partition $(\mu_1,\mu_2,\dotsc,\mu_{n-1},0)$ if
$\mu_n=0$, with the ``positive" partition
$(\mu_1,\dotsc,\mu_{n-1},\mu_n)\in Y_+$ if $\mu_n>0$, and with the
``negative" partition $(\mu_1,\dotsc,\mu_{n-1},-\mu_n)\in Y_-$ if
$\mu_n<0$.

Let $\mu$ and $\la$ be adjacent in $P$, and $x\in B(\Lab_1)$ such
that $\la-\mu=\wt(x)$. We think of this as walking from $\mu$ to
$\la$ by the step $x$. In terms of partitions, if $x=i$ for $1\le
i\le n-1$ then a cell is added to the $i$-th row. If $x=\ov{i}$ for
$1\le i\le n-1$ then a cell is removed from the $i$-th row. If $i=n$
then the above rules hold provided that $\la,\mu\in Y_+$. If
$\la,\mu\in Y_-$ then the roles of $n$ and $\ov{n}$ are reversed.

Let $B=B(\Lab_1)^{\otimes L}$. Let $b=b_L\dotsm b_1\in P(B)$ with
$b_j\in B(\Lab_1)$. In the usual way, $b$ can be regarded as a path
in the set of dominant weights: the $i$-th weight is given by the
weight of $b_i\dotsm b_1$. Alternatively $b$ describes a walk in $G$
from the empty partition to the element of $G$ corresponding to the
weight of $b$.

\begin{example}\label{ex:path} Let $n=4$. Consider $b=\ov{4}\ov{4}41321\in P(B)$
where $B=B(\Lab_1)^{\otimes 7}$. The element $b$ corresponds to the
walk in $G$ given by

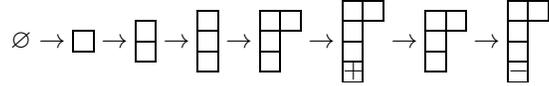
\begin{equation*}
\Yboxdim{8pt} \varnothing \rightarrow \yng(1) \rightarrow
\yng(1,1)\rightarrow \yng(1,1,1) \rightarrow \yng(2,1,1) \rightarrow
\young(\xx\xx,\xx,\xx,+) \rightarrow \yng(2,1,1) \rightarrow
\young(\xx\xx,\xx,\xx,-)
\end{equation*}
where the $+$ and $-$ markings on a partition indicate membership in
$Y_+$ and $Y_-$ respectively.
\end{example}

In the following proposition, for weights $\la,\mu\in P$, we write
$\mu\subset\la$ if the corresponding elements of $G$ are both in
$Y_+$ or both in $Y_-$ and the diagram of the partition associated
with $\mu$ is contained in that of $\la$.

\begin{prop}\label{pp:weight}
Suppose $b\in P(B^1\otimes B\otimes B^1,\la)$,
$\rh(b)\in P(B^1\otimes B,\alpha)$, $\lh(b)\in P(B\otimes
B^1,\beta)$ and $\rh(\lh(b))=\lh(\rh(b))\in P(B,\gamma)$. Then
$\alpha$ is uniquely determined by $\la$, $\beta$, and $\gamma$.
More precisely,
\begin{enumerate}
\item If $|\la|=|\gamma|+2$:
\begin{enumerate}
\item If the cells $\lambda/\beta$ and $\beta/\gamma$ are in
different rows and different columns, then $\alpha=\la-\{\beta/\gamma\}$.
\item If $\la/\beta$ and $\beta/\gamma$ are in the
same row or in the same column, then $\alpha=\beta$.
\end{enumerate}
\item If $|\la|=|\gamma|-2$:
\begin{enumerate}
\item If the cells $\beta/\la$ and $\gamma/\beta$ are in
different rows and different columns, then $\alpha=\la\cup\{\gamma/\beta\}$.
\item If $\beta/\lambda$ and $\gamma/\beta$ are in the
same row or the same column, then $\alpha=\beta$.
\end{enumerate}
\item If $|\la|=|\gamma|$ and $\la\not=\gamma$:
\begin{enumerate}
\item If $\la \supset \beta$ then $\alpha=\la\cup\{\gamma/\beta\}$.
\item If $\la \subset \beta$ then $\alpha=\la-\{\beta/\gamma\}$.
\end{enumerate}
\item If $\la=\gamma$:
\begin{enumerate}
\item If $\la \subset \beta$:
\begin{enumerate}
\item If $\beta/\la$ is in the first column of $\la$:
\begin{enumerate}
\item If $\beta/\la$ is in the $n$-th row, for $\beta\in Y_\pm$ let
$\alpha\in Y_{\mp}$ be the corresponding partition.
\item Otherwise let $\alpha=\beta$.
\end{enumerate}
\item Else $\alpha\subset \la$ and $\alpha$ is obtained from $\la$
by removing the corner cell in the column to the left of $\beta/\la$.
\end{enumerate}
\item If $\la\supset\beta$ then $\alpha$ is obtained from $\la$
by adjoining a cell to the column to the right of $\la/\beta$.
\end{enumerate}
\end{enumerate}
\end{prop}
\begin{proof} The rule for the weight $\alpha$ is given by van Leeuwen \cite[Rule
4.1.1]{vL:1998}: $\alpha$ is the unique dominant element in the Weyl
group orbit of the weight $\la+\gamma-\beta$. Using this rule the
proof is straightforward.
\end{proof}

\begin{remark} \label{rem:rhshape} The two operations
$\rh \circ \lh$ and $\lh \circ \rh$ define a pair of two-step walks
in the graph $G$ from $\la$ to $\gamma$, whose intermediate vertices
are $\beta$ and $\alpha$ respectively. If there is only one such
walk then $\alpha=\beta$; this occurs in cases (1b) and (2b). If
there are exactly two such walks then $\alpha$ is always chosen to
be the intermediate vertex not equal to $\beta$; this occurs in
cases (1a), (2a), (3a), and (3b). In the case that $\la=\gamma$
there may be many such walks; the proper choice of $\alpha$ given
$\beta$ is described in the proposition.
\end{remark}

\begin{example} Let $b$ be as in Example \ref{ex:path}. Then
$\lh(b)=\ov{4}41321$, $\rh(b)=2\ov{4}3121$, and
$\rh(\lh(b))=\lh(\rh(b))=\ov{4}3121$. Therefore $\la$ is the weight
$(2,1,1,-1)$ or the partition $(2,1,1,1)\in Y_-$, $\beta$ is the
weight and partition $(2,1,1,0)$, $\alpha$ is the weight
$(2,2,1,-1)$ and the partition $(2,2,1,1)\in Y_-$, and $\gamma$ is
the weight $(2,1,1,-1)$ and the partition $(2,1,1,1)\in Y_-$. Since
$\la=\gamma$ (as elements of $P$ or $G$) and $\beta \subset \la$ as
partitions, Case (4b) applies. The cell $\la/\beta$ is in the first
column; therefore $\alpha$ should be obtained from $\gamma$ by
adjoining a cell at the end of the second column, which agrees with
the example.
\end{example}

\begin{example} In $D_4^{(1)}$ let $b=4\ov{4}321\in
(B^{1,1})^{\otimes 5}$. Then $\lh(b)=\ov{4}321$, $\rh(b)=4321$,
$\rh(\lh(b))=\lh(\rh(b))=321$. Therefore $\la=(1,1,1,0)$, $\beta$ is
the weight $(1,1,1,-1)$ or the partition $(1,1,1,1)\in Y_-$,
$\gamma=(1,1,1,0)$, and $\alpha$ is the weight $(1,1,1,1)$ or the
partition $(1,1,1,1)\in Y_+$. This is case (4a1A).
\end{example}

\section{Rigged Configurations}\label{sec:RC}
In this section it is assumed that $\geh$ is nonexceptional and
simply-laced, that is, $\geh=A_n^{(1)}$ or $\geh=D_n^{(1)}$.

\subsection{Definition}\label{sec:rc def}
Let $B\in \CC$ for type $D_n^{(1)}$ and $B\in \CCA$ for type
$A_n^{(1)}$. Recall the notation in subsection \ref{ss:tensorcats},
where $L=(L_i^{(a)}\mid (a,i)\in \HH)$ is the multiplicity array of $B$.
The sequence of partitions $\nu=\{\nu^{(a)}\mid a\in \J \}$ is a
\textbf{$(L,\la)$-configuration} if
\begin{equation*}
\sum_{(a,i)\in\HH} i m_i^{(a)} \alpha_a = \sum_{(a,i)\in\HH} i
L_i^{(a)} \Lab_a- \la,
\end{equation*}
where $m_i^{(a)}$ is the number of parts of length $i$ in partition
$\nu^{(a)}$. A $(L,\la)$-configuration is \textbf{admissible} if
$p_i^{(a)}\ge 0$ for all $(a,i)\in\HH$, where $p_i^{(a)}$ is the
\textbf{vacancy number}
\begin{equation*}
p_i^{(a)}=\sum_{j\ge 1} \min(i,j) L_j^{(a)}
 - \sum_{b\in \J} (\alpha_a | \alpha_b) \sum_{j\ge 1}
 \min(i,j)m_j^{(b)}.
\end{equation*}
Here $(\cdot | \cdot )$ is the normalized invariant form on $P$
such that $(\alpha_i | \alpha_j)$ is the Cartan matrix. Let
$\Conf(L,\la)$ be the set of admissible $(L,\la)$-configurations. A
\textbf{rigged configuration} $(\nu,J)$ consists of a configuration
$\nu\in \Conf(L,\la)$ together with a double sequence of partitions
$J=\{J^{(a,i)}\mid (a,i)\in\HH \}$ such that the partition
$J^{(a,i)}$ is contained in a $m_i^{(a)}\times p_i^{(a)}$ rectangle.
The set of rigged configurations is denoted by $\RC(L,\la)$.

The partition $J^{(a,i)}$ is called \textbf{singular} if it has a
part of size $p_i^{(a)}$. The partition $J^{(a,i)}$ is called
\textbf{cosingular} if it has a part of size zero, or equivalently,
its complement in the rectangle of size $m_i^{(a)}\times p_i^{(a)}$
has a part of size $p_i^{(a)}$.

It is often useful to view a rigged configuration $(\nu,J)$ as a
sequence of partitions $\nu$ where the parts of size $i$ in
$\nu^{(a)}$ are labeled by the parts of $J^{(a,i)}$. The pair
$(i,x)$ where $i$ is a part of $\nu^{(a)}$ and $x$ is a part of
$J^{(a,i)}$ is called a \textbf{string} of the $a$-th rigged
partition $(\nu,J)^{(a)}$. The label $x$ is called a
\textbf{rigging} or \textbf{quantum number}. The corresponding
\textbf{coquantum number} is $p_i^{(a)}-x$.

\begin{example}\label{ex:rc}
Let $\geh=D_4^{(1)}$, $B=B^1\otimes B^2\otimes B^2\otimes B^3$ and
$\la=2\Lab_1$. Then the following three sequences of partitions
are admissible $(L,\la)$-configurations
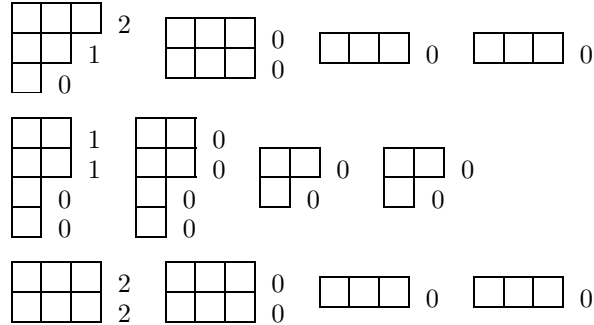
\begin{align*}
& \yngrc(3,2,2,1,1,0) \quad \yngrc(3,0,3,0) \quad \yngrc(3,0) \quad \yngrc(3,0)\\[2mm]
& \yngrc(2,1,2,1,1,0,1,0) \quad \yngrc(2,0,2,0,1,0,1,0)
  \quad \yngrc(2,0,1,0) \quad \yngrc(2,0,1,0)\\[2mm]
& \yngrc(3,2,3,2) \quad \yngrc(3,0,3,0) \quad \yngrc(3,0) \quad \yngrc(3,0)
\end{align*}
where the corresponding vacancy numbers are written next to each part. Hence, writing
the parts of $J^{(a,i)}$ next to the parts of size $i$ of partition $\nu^{(a)}$
the following would be a particular rigged configuration
\begin{equation*}
(\nu,J)\quad=\quad
\yngrc(3,0,2,1,1,0) \quad \yngrc(3,0,3,0) \quad \yngrc(3,0) \quad \yngrc(3,0).
\end{equation*}
\end{example}

\subsection{Quantum number complementation}
Let $\flip=\flip_L:\RC(L,\la)\rightarrow\RC(L,\la)$ be the
involution that preserves configurations and complements riggings
with respect to the vacancy numbers. More precisely, each partition
$J^{(a,i)}$ is replaced by the partition that is complementary to it
within the $m_i^{(a)}\times p_i^{(a)}$ rectangle.

\subsection{The RC reduction steps $\db$ and $\dt$}
\label{ss:deltas} Suppose $L_1^{(1)}>0$. Let $\Lh$ and $\Lr$
be obtained from $L$ by removing one tensor factor $B^1$. In
particular, if $B$ has $B^1$ as its left (resp. right) tensor
factor, then $\Lh$ (resp. $\Lr$) is the multiplicity array for $\lh(B)$ (resp.
$\rh(B)$). In \cite{OSS:2003b} a quantum number bijection
$\phib:P(B)\rightarrow \RC(L)$ was defined when $B$ is a tensor
power of $B^1$. The key step in the definition of $\phib$ is an
algorithm that defines a map $\db:\RC(L)\rightarrow\RC(\Lh)$. The
same algorithm defines such a map for the current case.

For $(\nu,J)\in\RC(L,\la)$, the algorithm produces a smaller rigged configuration
$\db(\nu,J)\in\RC(\Lh,\mu)$ for some $\mu\in\lm$ and an element
$\rk(\nu,J)\in B^1$ such that
\begin{equation} \label{eq:deltawt}
  \mu + \wt(\rk(\nu,J)) = \la.
\end{equation}

We recall the algorithm for $\db$ explicitly for type $A_n^{(1)}$
and $D_n^{(1)}$. Although we do not use them here, the explicit
algorithms exists for the other nonexceptional affine types and can
be found in~\cite{OSS:2002a}.

\subsubsection*{String selection for type $A_n^{(1)}$} Set $\ell^{(0)}=1$ and repeat the
following process
for $a=1,2,\ldots,n$ or until stopped. Find the smallest index $i\ge \ell^{(a-1)}$ such
that $J^{(a,i)}$ is singular. If no such $i$ exists, set $\rk(\nu,J)=a$ and stop.
Otherwise set $\ell^{(a)}=i$ and continue with $a+1$.

\subsubsection*{String selection for type $D_n^{(1)}$}
Set $\ell^{(0)}=1$ and repeat the following process for
$a=1,2,\ldots,n-2$ or until stopped. Find the smallest index $i\ge
\ell^{(a-1)}$ such that $J^{(a,i)}$ is singular. If no such $i$
exists, set $\rk(\nu,J)=a$ and stop. Otherwise set $\ell^{(a)}=i$
and continue with $a+1$. Set all yet undefined $\ell^{(a)}$ to
$\infty$.

If the process has not stopped at $a=n-2$,
find the minimal indices $i,j\ge \ell^{(n-2)}$ such that
$J^{(n-1,i)}$ and $J^{(n,j)}$ are singular. If neither $i$ nor
$j$ exist, set $\rk(\nu,J)=n-1$ and stop.
If $i$ exists, but not $j$, set $\ell^{(n-1)}=i$, $\rk(\nu,J)=n$ and stop.
If $j$ exists, but not $i$, set $\ell^{(n)}=j$, $\rk(\nu,J)=\overline{n}$
and stop. If both $i$ and $j$ exist, set $\ell^{(n-1)}=i$, $\ell^{(n)}=j$
and continue with $a=n-2$.

Now continue for $a=n-2,n-3,\ldots,1$ or until stopped.
Find the minimal index $i\ge \lb^{(a+1)}$ where $\lb^{(n-1)}
=\max(\ell^{(n-1)},\ell^{(n)})$ such that $J^{(a,i)}$ is singular
(if $i=\ell^{(a)}$ then there need to be two parts of size
$p_i^{(a)}$ in $J^{(a,i)}$).
If no such $i$ exists, set $\rk(\nu,J)=\overline{a+1}$ and stop.
If the process did not stop, set $\rk(\nu,J)=\overline{1}$.
Set all yet undefined $\ell^{(a)}$ and $\lb^{(a)}$ to $\infty$.

\subsubsection*{The new rigged configuration}
The rigged configuration $(\tilde{\nu},\tilde{J})=\db(\nu,J)$ is obtained by
removing a box from the selected strings and making the new strings singular
again. Explicitly (ignoring the statements about $\lb^{(a)}$ for type $A_n^{(1)}$)
\begin{equation*}
 m_i^{(a)}(\tilde{\nu})=m_i^{(a)}(\nu)+\begin{cases}
 1 & \text{if $i=\ell^{(a)}-1$}\\
 -1 & \text{if $i=\ell^{(a)}$}\\
 1 & \text{if $i=\lb^{(a)}-1$ and $1\le a\le n-2$}\\
 -1 & \text{if $i=\lb^{(a)}$ and $1\le a \le n-2$}\\
 0 & \text{otherwise.} \end{cases}
\end{equation*}
The partition $\tilde{J}^{(a,i)}$ is obtained from $J^{(a,i)}$ by removing
a part of size $p_i^{(a)}(\nu)$ for $i=\ell^{(a)}$ and $i=\lb^{(a)}$,
adding a part of size $p_i^{(a)}(\tilde{\nu})$ for $i=\ell^{(a)}-1$ and
$i=\lb^{(a)}-1$, and leaving it unchanged otherwise.

\begin{example}
For the rigged configuration $(\nu,J)$ of example~\ref{ex:rc}, we have
\begin{equation*}
\db(\nu,J)\quad=\quad
\yngrc(3,0,2,1) \quad \yngrc(2,0,2,0) \quad \yngrc(2,0) \quad \yngrc(2,0)
\end{equation*}
with $\rk(\nu,J)=\overline{2}$.
\end{example}

The next proposition was proved in~\cite{KSS:2002,OSS:2002a}.

\begin{prop} \label{pp:boxbij} The map
$\db:\RC(L,\la)\rightarrow \bigcup_{\mu\in \lm} \RC(\Lh,\mu)$ is
injective.
\end{prop}
Note that for simply-laced type, knowing $\la$ and $\mu$ uniquely
determines $\rk(\nu,J)$ by \eqref{eq:deltawt}.

We may define the inverse of $\db$. To this end, let
\begin{equation*}
\lp=\{ \mu\in\Pb^+ \mid \text{$B(\mu)$ occurs in $B^1\times B(\la)$} \}.
\end{equation*}
Denote by $\RCt(L,\la)$ the subset of $\RC(L,\la)\times B^1$ given by
$((\nu,J),b)$ such that $\la+\wt(b)\in \Pb^+$. By abuse of notation
define
\begin{equation*}
\db^{-1}:\RCt(L,\la)\to \bigcup_{\beta\in \lp}\RC(\Lp,\beta)
\end{equation*}
by the following algorithm, where $\Lp$ is obtained from $L$ by replacing
$L_1^{(1)}$ by $L_1^{(1)}+1$.

\subsubsection*{String selection for type $A_n^{(1)}$} In this case $\wt(b)=\epsilon_r$
for some $1\le r\le n+1$, where $\epsilon_r$ is the $r$-th canonical unit vector in
$\Z^{n+1}$. Set $s^{(r)}=\infty$ and repeat the following
process for $a=r-1,r-2,\ldots,1$. Find the largest index $i\le s^{(a+1)}$ such
that $J^{(a,i)}$ is singular and set $s^{(a)}=i$; if no such $i$ exists set $s^{(a)}=0$.
Set all undefined $s^{(a)}$ to infinity.

\subsubsection*{String selection for type $D_n^{(1)}$}
In this case $\wt(b)=\epsilon_r$ or $\wt(b)=-\epsilon_r$ for $1\le r\le n$, where
$\epsilon_r$ is the $r$-th canonical unit vector in $\Z^{n}$.
In the first case proceed exactly as for type $A_n^{(1)}$. Throughout the
whole algorithm, if an index $i$ does not exist, set $i=0$.

If $\wt(b)=-\epsilon_n$, find the largest index $i$ such that $J^{(n,i)}$ is singular
and set $s^{(n)}=i$. Find the largest index $i\le s^{(n)}$ such that $J^{(n-2,i)}$ is
singular and set $s^{(n-2)}=i$. Then proceed as in type $A_n^{(1)}$.

If $\wt(b)=-\epsilon_{n-1}$, find the largest indices $i$ and $j$ such that
$J^{(n-1,i)}$ and $J^{(n,j)}$ are singular and set $s^{(n-1)}=i$ and $s^{(n)}=j$.
Then find the largest index $i\le \min\{s^{(n-1)},s^{(n)}\}$ such that $J^{(n-2,i)}$
is singular and set $s^{(n-2)}=i$. After this proceed as in type $A_n^{(1)}$.

Finally, if $\wt(b)=-\epsilon_r$ for $1\le r\le n-2$, set $\sbar^{(r-1)}=\infty$ and
proceed for $a=r,r+1,\ldots,n-2$ as follows. Find the largest index $i\le \sbar^{(a-1)}$
such that $J^{(a,i)}$ is singular and set $\sbar^{(a)}=i$. Then find the largest indices
$i\le \sbar^{(n-2)}$ and $j\le \sbar^{(n-2)}$ such that
$J^{(n-1,i)}$ and $J^{(n,j)}$ are singular and set $s^{(n-1)}=i$ and $s^{(n)}=j$.
After this proceed as for the case $\wt(b)=-\epsilon_{n-1}$.

Set all yet undefined $s^{(a)}$ and $\sbar^{(a)}$ to $\infty$.

\subsubsection*{The new rigged configuration}
The rigged configuration $(\tilde{\nu},\tilde{J})=\db^{-1}(\nu,J)$ is obtained by
adding a box to the selected strings and making the new strings singular
again.

Define $\dt:\RC(L)\rightarrow\RC(\Lh)$ by
$\flip_{\Lh}\circ\db\circ\flip_L$. Alternatively, $\dt$ is defined
by a coquantum number version of the map $\db$. Instead of selecting
singular strings it selects cosingular strings and keeps coquantum
numbers constant for unselected strings. It also produces an element
$\rkt(\nu,J)\in B^1$. If $(\nu,J)\in\RC(L,\la)$ and $\dt(\nu,J)\in
\RC(\Lh,\mu)$ then
\begin{equation*}
  \mu + \wt(\rkt(\nu,J)) = \la.
\end{equation*}

\subsection{Splitting on RCs}
 Let $s\ge2$. Suppose $B$ contains a distinguished tensor factor
$B^{r,s}$, which is the case when we consider the maps $\ls$ and
$\rs$. Let $L$ be the multiplicity array of $B$ and $\Ls$ that which
is obtained from $L$ by replacing $B^{r,s}$ by $B^{r,1}$ and
$B^{r,s-1}$.

\begin{prop} \label{pp:ljdef} Let $L$ be such that $L_s^{(r)}\ge1$ for a particular
$(r,s)\in\HH$ with $s\ge2$ and let $\Ls$ be defined with respect to
$(r,s)$. Then $\Conf(L,\la)\subset \Conf(\Ls,\la)$. Under this
inclusion map, the vacancy number $p_i^{(a)}$ for $\nu$ increases by
$\delta_{a,r} \chi(i<s)$ where $\chi(P)=1$ if $P$ is true and
$\chi(P)=0$ otherwise. Hence there are well-defined injective maps
$\rcls,\rcrs:\RC(L)\rightarrow \RC(\Ls)$ given by:
\begin{enumerate}
\item $\rcls(\nu,J)=(\nu,J)$.
\item $\rcrs(\nu,J)=(\nu,J')$ where $J'$ is obtained from $J$ by
adding $1$ to the rigging of each string in $(\nu,J)^{(r)}$ of
length strictly less than $s$.
\end{enumerate}
In particular, $\rcls$ preserves quantum numbers, $\rcrs$ preserves
coquantum numbers, and
\begin{equation}\label{eq:rj-lj}
  \rcrs = \flip_{\Ls} \circ \rcls \circ \flip_L.
\end{equation}

\end{prop}
\begin{proof} Immediate from the definitions.
\end{proof}

\subsection{Box-splitting for RCs} Suppose $r\ge2$ and $B\in\CCA$ has a
distinguished tensor factor $B^{r,1}$. Let $L$ be the multiplicity
array for $B$ and $\Lb$ that for the crystal obtained from $B$ by
replacing $B^{r,1}$ by $B^{1,1}$ and $B^{r-1,1}$.

\begin{prop} \label{pp:rcbs} Let $L$ be such that $L_1^{(r)}\ge1$
for some $r\ge2$. Let $\Lb$ be defined with respect to $r$. Then
there are injections $\rclb,\rcrb:\RC(L,\la)\rightarrow
\RC(\Lb,\la)$ defined by adding singular (resp. cosingular) strings
of length $1$ to $(\nu,J)^{(a)}$ for $1\le a < r$. Moreover the
vacancy numbers stay the same.
\end{prop}

\section{Fermionic formula $M$}\label{s:fermionic}
In this section we state the fermionic formula $M$ associated with
rigged configurations for simply-laced algebras as introduced in~\cite{HKOTY:1999}
and virtual fermionic formulas for nonsimply-laced algebras
(see \cite{OSS:2003a,OSS:2003b}).

\subsection{Fermionic formula $M$}
Let $(q)_m=(1-q)(1-q^2)\cdots (1-q^m)$ and define the $q$-binomial coefficient
for $m,p\in \Z_{\ge 0}$ as
\begin{equation*}
\qbin{m+p}{m}=\frac{(q)_{m+p}}{(q)_m(q)_p}.
\end{equation*}
The fermionic formula for types $A_n^{(1)}$ and $D_n^{(1)}$ is given by~\cite{HKOTY:1999}
\begin{equation}\label{eq:fermionic}
M_{L,\la}(q)=\sum_{\nu\in \Conf(L,\la)} q^{\cc(\nu)}
\prod_{(a,i)\in\HH} \qbin{m_i^{(a)}+p_i^{(a)}}{m_i^{(a)}}
\end{equation}
with $m_i^{(a)}$, $p_i^{(a)}$ and $\Conf(L,\la)$ as in section~\ref{sec:rc def}
and
\begin{equation*}
\cc(\nu)=\frac{1}{2} \sum_{a,b\in \J} \sum_{j,k\ge 1} (\alpha_a \mid
\alpha_b) \min(j,k) m_j^{(a)} m_k^{(b)}.
\end{equation*}
The fermionic formula \eqref{eq:fermionic} can be restated solely in terms of
rigged configurations. To this end recall that the $q$-binomial coefficient
$\qbin{m+p}{m}$ is the generating function of partitions in a box of width $p$
and height $m$. Hence
\begin{equation}\label{eq:M rc}
M_{L,\la}(q)=\sum_{(\nu,J)\in \RC(L,\la)} q^{\cc(\nu,J)},
\end{equation}
where $\cc(\nu,J)=\cc(\nu)+\sum_{(a,i)\in \HH} |J^{(a,i)}|$.

\subsection{Virtual fermionic formula}
Fermionic formulae for nonsimply-laced algebras were defined
in~\cite[Section 4]{HKOTT:2001}. For $A_{2n}^{(2)\dagger}$ it was
defined in \cite{OSS:2003a}. Here we recall virtual rigged
configurations in analogy to virtual crystals as defined
in~\cite{OSS:2003b}.

\begin{definition} \label{def:VRC}
Let $X$ and $Y$ be as in section~\ref{ss:VX}, and $\la$, $B$ and $L$
as in section \ref{sec:rc def} for type $X$. Let $\emb:B\to\Vh$ be
the corresponding virtual $Y$-crystal and $\Lhat$ the multiplicity array corresponding
to $\Vh$. For
$X\not\in\{A_{2n}^{(2)},A_{2n}^{(2)\dagger}\}$, $\RCv(L,\la)$ is the
set of elements $(\nh,\Jh)\in \RC(\Lhat,\emb(\la))$ such that:
\begin{enumerate}
\item For all $i\in \Z_{>0}$,
$\mh_i^{(a)}=\mh_i^{(b)}$ and $\Jh^{(a,i)}=\Jh^{(b,i)}$ if $a$ and $b$ are in the
same $\aut$-orbit in $I^Y$.
\item For all $i\in \Z_{>0}$, $a\in \J^X$, and $b\in \bij(a)\subset \J^Y$, we have
$\mh_j^{(b)}=0$ if $j \not\in \mult_a \Z$ and the parts of $\Jh^{(b,i)}$ are multiples
of $\mult_a$.
\end{enumerate}
For $X=A_{2n}^{(2)}$ the following changes must be made:
\begin{enumerate}
\item[(A2)] $\mh_j^{(n)}$ may be positive for any $j\ge1$.
\end{enumerate}
For $X=A_{2n}^{(2)\dagger}$ one makes the exception (A2) and the
additional condition that
\begin{enumerate}
\item[(A2D)] The parts of $\Jh^{(n,i)}$ must have the same parity
as $i$.
\end{enumerate}
\end{definition}

\begin{theorem}\cite[Theorem 4.2]{OSS:2003b} \label{th:M=VM}
There is a bijection $\Psi:\RC(L,\la)\rightarrow \RCv(L,\la)$
sending $(\nu,J)\mapsto(\nh,\Jh)$ given as follows. For all $a\in
\J^X$, $b\in\bij(a)\subset \J^Y$, and $i\in\Z_{>0}$,
\begin{align*}
\mh_{\mult_a i}^{(b)} &= m_i^{(a)} \\
\Jh^{(b,\mult_a i)}&=\mult_a J^{(a,i)},
\end{align*}
except when $X=A_{2n}^{(2)}$ or $X=A_{2n}^{(2)\dagger}$ and $a=n$,
in which case
\begin{equation*}
\begin{split}
  \mh_i^{(n)} &= m_i^{(n)} \\
  \Jh^{(n,i)} &= 2 J^{(n,i)}.
\end{split}
\end{equation*}
The cocharge changes by $\cc(\nh,\Jh) = \mult_0 \,\cc(\nu,J)$.
\end{theorem}
Defining the virtual fermionic formula as
\begin{equation*}
VM_{L,\la}(q)=\sum_{(\nh,\Jh)\in\RCv(L,\la)} q^{\cc(\nh,\Jh)/\mult_0}
\end{equation*}
we obtain as a corollary:
\begin{corollary} \label{cor:M=VM} ($M=VM$)
$M_{L,\la}(q)=VM_{L,\la}(q)$.
\end{corollary}

\section{Bijection} \label{sec:bij}

\subsection{Quantum number bijection}
The following result defines the bijection from paths to rigged
configurations. It is valid for both $B\in \CC$ and $B\in \CCA$.

\begin{prop} \label{pp:phi} There exists a unique family of bijections
$\phib:P(B,\la)\rightarrow \RC(L,\la)$ such that the empty path maps
to the empty rigged configuration, and:
\begin{enumerate}
\item Suppose $B=B^1 \otimes B'$. Let $\lh(B)=B'$ with multiplicity array $\Lh$. Then the diagram
\begin{equation} \label{eq:hatdiag}
\begin{CD}
P(B,\la) @>{\phib}>> \RC(L,\la) \\
@V{\lh}VV @VV{\db}V \\
\displaystyle{\bigcup_{\mu\in\lm} P(\lh(B),\mu)} @>>{\phib}> \displaystyle{\bigcup_{\mu\in\lm}
\RC(\Lh,\mu)}
\end{CD}
\end{equation}
commutes.
\item Suppose $B=B^{r,s} \otimes B'$ with $s\ge 2$ (and $r=1$ for $\CC$).
Let $\ls(B)=B^{r,1}\otimes B^{r,s-1}\otimes B'$ with multiplicity
array $\Ls$. Then the following diagram commutes:
\begin{equation} \label{eq:splitdiag}
\begin{CD}
P(B,\la) @>{\phib}>> \RC(L,\la) \\
@V{\ls}VV @VV{\rcls}V \\
P(\ls(B),\la) @>>{\phib}> \RC(\Ls,\la)
\end{CD}
\end{equation}
\item For $\CCA$, suppose $B=B^{r,1} \otimes B'$ with $r\ge2$. Let $\llb(B)=B^{1,1}\otimes B^{r-1,1}\otimes B'$
and $\Lb$ its multiplicity array. Then the following diagram
commutes:
\begin{equation}\label{eq:boxdiag}
\begin{CD}
P(B,\la) @>{\phib}>> \RC(L,\la) \\
@V{\llb}VV @VV{\rclb}V \\
P(\llb(B),\la) @>>{\phib}> \RC(\Lb,\la)
\end{CD}
\end{equation}
\end{enumerate}
\end{prop}

For type $A_n^{(1)}$ the existence of $\phib$ was proven
in~\cite{KSS:2002}. The proof in case (1) for other nonexceptional
types is essentially done in~\cite{OSS:2002a}. It remains to prove
case (2) for type $D_n^{(1)}$.

\begin{lemma} \label{lem:resls}
Let $B=B^s\otimes B'$ with $s\ge 2$. For type $D_n^{(1)}$,
the map $\phib:P(\ls(B),\la)\to \RC(\Ls,\la)$ restricts to
a bijection $\phib:\ls(P(B,\la))\to \rcls(\RC(L,\la))$.
\end{lemma}

\begin{proof}
Let $b=x\otimes b_2 \otimes b'\in B^1\otimes B^{s-1}\otimes B'$ and
$\ls(b_2)=y\otimes b_3\in B^1\otimes B^{s-2}$. Then $b\in
\Image(\ls)$ if and only if $x\le y$. (Note that this implies in
particular that $n$ and $\bar{n}$ cannot appear in the same one-row
crystal element).

By Proposition~\ref{pp:ljdef}, $(\nu,J)\in \RC(\Ls,\la)$ is in
the image of $\rcls$ if and only if $(\nu,J)^{(1)}$ has
no singular strings of length smaller than $s$.

Let us first show that if $b\in\Image(\ls)$ then
$\phib(b)\in \Image(\rcls)$. Hence assume that $b=x\otimes b_2\otimes b'$
with $x\le y$ with $y$ as defined above. By induction $(\nu',J')=
\phib(y\otimes b_3\otimes b')$ has no singular strings in the first rigged
partition of length smaller than $s-1$. Denote the lengths of the strings
selected by $\db$ associated with the letter $y$ by $\ell_y^{(k)}$ and
$\lb_y^{(k)}$. Then in particular $\ell_y^{(1)}\ge s-1$. ``Unsplitting'' yields
on the paths side $b_2\otimes b'$ and on the rigged configuration side $(\nu',J')$
with a change in the vacancy numbers by $-\delta_{a,1}\chi(i<s-1)$. Since $x\le y$
it follows that $\ell_x^{(k)}>\ell_y^{(k)}$ and $\lb_x^{(k)}>\lb_y^{(k)}$, where
$\ell_x^{(k)}$ and $\lb_x^{(k)}$ are the lengths of the strings selected by $\db$
associated with $x$. This shows in particular that $\ell_x^{(1)}\ge s$, and
from the change in vacancy numbers from $\phib(b_2\otimes b')$ to
$\phib(x\otimes b_2\otimes b')$ it follows that there are no singular
strings in the first rigged partition of $\phib(x\otimes b_2\otimes b')$
of length smaller than $s$.

Conversely, assume that $(\nu,J)\in \RC(\Ls,\la)$ is in the image if $\rcls$.
We need to show that then $b=\phib^{-1}(\nu,J)$ has the property that
$x\le y$ in the above notation. Call the strings selected
by $\db$ in $(\nu,J)$ $\ell_x^{(k)}$ and $\lb_x^{(k)}$. By assumption
$(\nu,J)^{(1)}$ has no singular string of length smaller than $s$.
Hence $\ell_x^{(1)}\ge s$. By the definition of $\rcls$, we have that
the first rigged partition of $(\nu',J')=\rcls\circ \db(\nu,J)$ has no
singular strings of length smaller than $s-1$. Hence $s-1\le \ell_y^{(1)}<\ell_x^{(1)}$,
where $\ell_y^{(k)}$ and $\lb_y^{(k)}$ are the lengths of the strings selected
by $\db$ on $(\nu',J')$. The algorithm of $\db$ implies that $\ell_y^{(k)}<\ell_x^{(k)}$
and $\lb_y^{(k)}<\lb_x^{(k)}$, so that $x\le y$ as desired.
\end{proof}

\subsection{Coquantum number bijection}
Let $\phit =\flip\circ\phib$; it can be characterized as follows.

\begin{prop} There exists a unique family of bijections
$\phit:P(B,\la)\rightarrow \RC(L,\la)$ with the same properties as in
Proposition~\ref{pp:phi} except that $\db$, $\rcls$ and $\rclb$
are replaced by $\dt$, $\rcrs$ and $\rcrb$ in \eqref{eq:hatdiag},
\eqref{eq:splitdiag} and \eqref{eq:boxdiag}, respectively.
\end{prop}

\subsection{Commutations of the basic steps}
We record the commutations among the basic steps of the path-RC
bijection. Here $r=1$ for $\CC$.

\begin{theorem} \label{th:RCcomm}
\mbox{}
\begin{enumerate}
\item \label{pt:dbdt}
$[\db,\dt]=0$.
\item \label{pt:jdelta} $[\rcrs,\db]=0$ and $[\rcls,\dt]=0$.
\item \label{pt:ljrj} $[\rcls,\rcrs]=0$.
\item \label{pt:bsdelta} $[\rcrb,\db]=0$ and $[\rclb,\dt]=0$.
\item \label{pt:bslj} $[\rcrb,\rcls]=0$ and $[\rclb,\rcrs]=0$.
\end{enumerate}
\end{theorem}

The proof of part \ref{pt:dbdt} for type $A_n^{(1)}$ is given
in~\cite[Appendix A]{KSS:2002}. The proof of part~\ref{pt:dbdt}
for type $D_n^{(1)}$ is quite technical and follows similar arguments
as~\cite[Appendix A]{KSS:2002} (see also~\cite[Appendix C]{Sch:2004}). Details are
available upon request. Parts \ref{pt:jdelta} and \ref{pt:ljrj}
follow easily from the definitions. Parts \ref{pt:bsdelta} through
\ref{pt:bslj} only apply for $\CCA$ and follow from \cite{KSS:2002}.

For type $D_n^{(1)}$, there is an analogue of
Proposition~\ref{pp:weight} for the commutation of $\db$ and $\dt$.
Let $(\nu,J)\in \RC(L,\la)$, $\dt(\nu,J)\in \RC(\rh(L),\alpha)$,
$\db(\nu,J)\in \RC(\lh(L),\beta)$ and
$\dt(\db(\nu,J))=\db(\dt(\nu,J))\in \RC(\lh(\rh(L)),\gamma)$. Then
$\alpha$ is uniquely determined by $\la$, $\beta$, and $\gamma$.
\begin{prop}\label{pp:weight rc}
For $\la$, $\alpha$, $\beta$, and $\gamma$ as above the statements
of Proposition~\ref{pp:weight} hold.
\end{prop}
The proof is an easy consequence of the commutation $[\db,\dt]=0$ and is
available upon request.

\subsection{The bijection and the various operations}
\begin{theorem} \label{th:bijops}
Under the family of bijections $\phib$ the following operations
correspond:
\begin{enumerate}
\item \label{pt:ls} $\ls$ with $\rcls$.
\item \label{pt:lh} $\lh$ with $\db$.
\item \label{pt:rs} $\rs$ with $\rcrs$.
\item \label{pt:rh} $\rh$ with $\dt$.
\item \label{pt:*} $*$ with $\flip$.
\item \label{pt:R} $R$ with the identity.
\item \label{pt:bs} $\llb$ with $\rclb$ and $\rb$ with $\rcrb$.
\end{enumerate}
\end{theorem}

\begin{example}
To illustrate point~\ref{pt:rh} of the above Theorem, take
\begin{equation*}
b=\begin{array}{|c|}\hline \overline{3}\\ \hline \end{array} \otimes
\young(23) \otimes \young(12) \otimes \young(1)
\end{equation*}
of type $D_4^{(1)}$. Then
\begin{equation*}
\rh(b)=\young(3) \otimes \young(22) \otimes \young(11)
\end{equation*}
and
\begin{align*}
\phib(b)&= \quad \yngrc(2,0,1,0,1,0) \quad \yngrc(1,1,1,0) \quad
\yngrc(1,0)
\quad \yngrc(1,0)\\[2mm]
\phib(\rh(b))=\dt(\phib(b))&= \quad \yngrc(2,0,1,0) \quad \yngrc(1,0) \quad \es \quad
\qquad \es.
\end{align*}
\end{example}

\begin{proof} Everything is proved for $\CCA$ in \cite{KSS:2002},
including part \ref{pt:bs}, which only applies in that case. We
assume that $B\in\CC$ for type $D_n^{(1)}$. Parts \ref{pt:ls} and
\ref{pt:lh} hold by Proposition \ref{pp:phi}. We prove parts
\ref{pt:rs}, \ref{pt:rh}, and \ref{pt:*} simultaneously by
induction. The induction is based first on the quantity $\sum_i s_i$
for the crystal $\bigotimes_i B^{s_i}$, and then by decreasing
induction on the number of tensor factors.

\Yboxdim{5pt}
\Yinterspace{1pt}

Consider part \ref{pt:rs}. Suppose first that $B=B^s$ for some $s\ge
2$. Then $P(B)$ has only one element $1^s$. It is easy to show that
$\phib(1^s)$ is the empty RC and that \ref{pt:rs} holds. Suppose
next that $B=B^1 \otimes B' \otimes B^s$. Consider the diagram
\Yboxdim{5pt} \Yinterspace{1pt}
\begin{equation}
\xymatrix{
P(B) \ar[rrr]^{\rs} \ar[dr]_{\phib} \ar[ddd]_{\lh} & & & P(\rs(B)) \ar[dl]^{\phib} \ar[ddd]^{\lh} \\
& \RC(L) \ar[r]^{\rcrs} \ar[d]_{\db} & \RC(\Lrs) \ar[d]^{\db} & \\
& \RC(\Lh) \ar[r]_{\rcrs} & \RC(\rs(\Lh)) & \\
P(\lh(B)) \ar[ur]^{\phib} \ar[rrr]_{\rs} & & & P(\rs(\lh(B)))
\ar[ul]_{\phib} }
\end{equation}
Here $L$, $\Lrs$, $\Lh$, $\rs(\Lh)$ are the multiplicities arrays
corresponding to $B$, $\rs(B)$, $\lh(B)$, $\rs(\lh(B))$, respectively.
We shall view such a diagram as a cube in which the small square is
in the background. The left and right faces commute by Proposition
\ref{pp:phi}. The front and back faces commute by Proposition
\ref{pp:cryscomm} part \ref{pt:lhrs} and Theorem \ref{th:RCcomm}
part \ref{pt:jdelta} respectively. The bottom face commutes by
induction. It follows that the top face ``commutes up to $\db$",
that is, $\db \circ \rcrs \circ \phib = \db \circ \phib \circ \rs$.
But all maps in the top face preserve the highest weight. By
Proposition \ref{pp:boxbij} it follows that the top face commutes.

The remaining case is $B=B^{s'} \otimes B' \otimes B^s$ for $s,s'\ge
2$. Consider the diagram below, where $\rs(\ls(L))$
is obtained from $\Ls$ by splitting a $B^{s'}$ into $B^{s'-1}$ and $B^1$.

\Yboxdim{5pt}
\Yinterspace{1pt}
\begin{equation}
\xymatrix{
P(B) \ar[rrr]^{\rs} \ar[dr]_{\phib} \ar[ddd]_{\ls} & & & P(\rs(B)) \ar[dl]^{\phib} \ar[ddd]^{\ls} \\
& \RC(L) \ar[r]^{\rcrs} \ar[d]_{\rcls} & \RC(\rs(L)) \ar[d]^{\rcls} & \\
& \RC(\ls(L)) \ar[r]_{\rcrs} & \RC(\rs(\ls(L)) & \\
P(\ls(B)) \ar[ur]^{\phib} \ar[rrr]_{\rs} & & & P(\rs(\ls(B)))
\ar[ul]_{\phib} }
\end{equation}
The left and right faces commute by Proposition \ref{pp:phi}. The
front and back faces commute by Proposition \ref{pp:cryscomm} part
\ref{pt:lsrs} and Theorem \ref{th:RCcomm} part \ref{pt:ljrj}
respectively. The bottom face commutes by induction. Since $\rcls$
is injective, it follows that the top face commutes. This finishes
the proof of part \ref{pt:rs}.

We now prove part \ref{pt:rh}. The proof is trivial for the base
case $B=B^1$. Suppose next that $B=B^1 \otimes B' \otimes B^1$.
\Yboxdim{5pt} \Yinterspace{1pt}
\begin{equation}
\xymatrix{
P(B) \ar[rrr]^{\rh} \ar[dr]_{\phib} \ar[ddd]_{\lh} & & & P(\rh(B)) \ar[dl]^{\phib} \ar[ddd]^{\lh} \\
& \RC(L) \ar[r]^{\dt} \ar[d]_{\db} & \RC(\Lh) \ar[d]^{\db} & \\
& \RC(\Lh) \ar[r]_{\dt} & \RC(\rh(\Lh)) & \\
P(\lh(B)) \ar[ur]^{\phib} \ar[rrr]_{\rh} & & & P(\rh(\lh(B)))
\ar[ul]_{\phib} }
\end{equation}
The left and right faces commute by Proposition \ref{pp:phi}. The
front and back faces commute by Proposition \ref{pp:cryscomm} part
\ref{pt:lhrh} and Theorem \ref{th:RCcomm} part \ref{pt:dbdt}
respectively. The bottom face commutes by induction. Thus the top
face commutes up to $\db$. By Proposition \ref{pp:boxbij} it
suffices to show that both ways around the top face, result in
elements with the same highest weight.
This follows from Propositions~\ref{pp:weight} and~\ref{pp:weight
rc}.

The remaining case is $B=B^s \otimes B' \otimes B^1$. Consider the
diagram \Yboxdim{5pt} \Yinterspace{1pt}
\begin{equation}
\xymatrix{
P(B) \ar[rrr]^{\rh} \ar[dr]_{\phib} \ar[ddd]_{\ls} & & & P(\rh(B)) \ar[dl]^{\phib} \ar[ddd]^{\ls} \\
& \RC(L) \ar[r]^{\dt} \ar[d]_{\rcls} & \RC(\Lr) \ar[d]^{\rcls} & \\
& \RC(\Ls) \ar[r]_{\dt} & \RC(\rh(\Ls)) & \\
P(\ls(B)) \ar[ur]^{\phib} \ar[rrr]_{\rh} & & & P(\rh(\ls(B)))
\ar[ul]_{\phib} }
\end{equation}
The left and right faces commute by Proposition \ref{pp:phi}. The
front and back faces commute by Proposition \ref{pp:cryscomm} part
\ref{pt:rhls} and Theorem \ref{th:RCcomm} part \ref{pt:jdelta}
respectively. The bottom face commutes by induction. Since $\rcls$
is injective it follows that the top face commutes. This proves part
\ref{pt:rh}.

For part \ref{pt:*} the proof of the base case $B=B^s$ is trivial.
Suppose next that $B=B^1 \otimes B' \otimes B^1$. Consider the
diagram \Yboxdim{5pt} \Yinterspace{1pt}
\begin{equation}
\xymatrix{
P(B) \ar[rrr]^{*} \ar[dr]_{\phib} \ar[ddd]_{\rh} & & & P(B^*) \ar[dl]^{\phib} \ar[ddd]^{\lh} \\
& \RC(L) \ar[r]^{\flip} \ar[d]_{\dt} & \RC(L) \ar[d]^{\db} & \\
& \RC(\Lr) \ar[r]_{\flip} & \RC(\Lr) & \\
P(\rh(B)) \ar[ur]^{\phib} \ar[rrr]_{*} & & & P(\rh(B)^*)
\ar[ul]_{\phib} }
\end{equation}
The right face commutes by Proposition~\ref{pp:phi}. The left
commutes by part \ref{pt:rh} which was just proved above. The back
face commutes by the definition of $\dt$. The commutation of the
front face is given by Proposition \ref{pp:cryscomm} part
\ref{pt:hat*}. The bottom face commutes by induction. It follows
that the top face commutes up to $\db$. Again it suffices to show
that both ways around the top face produce elements of the same
highest weight. But this holds since $\phib$, $\flip$, and $*$
preserve the highest weight. Here we are using the fact that for
$\la\in \Pb^+$, $V_\la^*\cong V_\la$.

Next let $B=B' \otimes B^s$ with $s\ge2$. \Yboxdim{5pt}
\Yinterspace{1pt}
\begin{equation}
\xymatrix{
P(B) \ar[rrr]^{*} \ar[dr]_{\phib} \ar[ddd]_{\rs} & & & P(B^*) \ar[dl]^{\phib} \ar[ddd]^{\ls} \\
& \RC(L) \ar[r]^{\flip} \ar[d]_{\rcrs} & \RC(L) \ar[d]^{\rcls} & \\
& \RC(\Ls) \ar[r]_{\flip} & \RC(\Ls) & \\
P(\rs(B)) \ar[ur]^{\phib} \ar[rrr]_{*} & & & P(\rs(B)^*)
\ar[ul]_{\phib} }
\end{equation}
The right face commutes by Proposition \ref{pp:phi}. The left face
commutes by part \ref{pt:rs} which was proved above. The back face
commutes by the definition of $\rcrs$. The commutation of the
front face is given by the definition of $\ls$ in
\eqref{eq:splitstar}. The bottom face commutes by induction. Since
$\rcls$ is injective, the top face commutes.

For $B=B^s \otimes B'$ with $s\ge2$ the proof is similar to the
previous case.

This concludes the proof of part \ref{pt:*}.

For the proof of part \ref{pt:R}, let $B=B_k\otimes B_{k-1}\otimes \cdots \otimes B_1$.
We may assume that $R=R_j$ is the
R-matrix being applied at tensor positions $j$ and $j+1$ (from the
right). By induction we may assume that $j=k-1$, that is, $R$ acts
at the leftmost two tensor positions. By part \ref{pt:*} and
Proposition \ref{pp:R*} we may assume that $j=1$. Again by induction
we may assume that $k=2$. Let $B=B^t \otimes B^s$ (of type
$D_n^{(1)}$). By Lemma \ref{lem:D2P} $B$ is multiplicity-free as a
$U_q(D_n)$-crystal. Since $R$ preserves weights it follows that
$R(v_{p,q}^{t,s})=v_{p,q}^{s,t}$. A direct computation shows that
$\phib(v_{p,q}^{t,s})=\phib(v_{p,q}^{s,t})$.
\end{proof}

\subsection{$X=M$ for types $A_n^{(1)}$ and $D_n^{(1)}$}
In this subsection we will show that $X_{B,\la}=M_{L,\la}$ for
$B\in\CCA$ for type $A_n^{(1)}$ and $B\in\CC$ for type $D_n^{(1)}$.
By Proposition~\ref{pp:phi} there is a bijection between the sets
$P(B,\la)$ and $\RC(L,\la)$. Hence it remains to show that the statistics
is preserved.

\begin{theorem}\label{th:statistics}
Let $B\in\CCA$ be a crystal of type $A_n^{(1)}$ or $B\in\CC$ a crystal of
type $D_n^{(1)}$ and $\la$ a dominant integral weight.
The coquantum number bijection $\phit$ preserves the statistics, that is
$D_B(b)=\cc(\phit(b))$ for all $b\in P(B,\la)$.
\end{theorem}

\begin{proof}
For type $A_n^{(1)}$ the theorem follows
from~\cite[Theorem 9.1]{KSS:2002}. Hence assume that $B\in\CC$ of
type $D_n^{(1)}$.
By Theorem~\ref{th:bijops} part \ref{pt:rs} and
equations~\eqref{eq:splitstar} and~\eqref{eq:rj-lj} the maps $\rs$
and $\rcls$ correspond under $\phit$. By Theorem~\ref{th:splitrow}
we have $D(\rs(b))=D(b)$. Similarly, it follows immediately from the
definition of $\rcls$ in Proposition~\ref{pp:ljdef} that
$\cc(\rcls(\nu,J))=\cc(\nu,J)$. The maps $R$ and the identity also
correspond under $\phit$ by Theorem~\ref{th:bijops} part \ref{pt:R},
and neither of them changes the statistics.

There exists a sequence $\mathcal{S}_P$ of maps $\rs$ and $R$ which
transforms a path $b\in P(B,\la)$ into a path of single boxes. By
Theorem~\ref{th:bijops} there exists a corresponding sequence
$\mathcal{S}_{\RC}$ of maps $\rcls$ and the identity. Since neither
of these maps changes the statistics it follows that
\begin{equation*}
D(\mathcal{S}_P(b))=\cc(\mathcal{S}_{\RC}(\phit(b))) \qquad \text{implies that}
\qquad D(b)=\cc(\phit(b)).
\end{equation*}
The theorem for the case $B=(B^{1,1})^{\otimes N}$ has already been proven
in~\cite{OSS:2002a}.
\end{proof}

\begin{corollary}\label{cor:X=M AD}
For $B\in\CCA$ of type $A_n^{(1)}$ or $B\in\CC$ of type $D_n^{(1)}$,
$L$ the corresponding multiplicity array and $\la$ a dominant integral,
we have
\begin{equation*}
X_{B,\la}(q)=M_{L,\la}(q).
\end{equation*}
\end{corollary}
\begin{proof}
This follows from Theorem~\ref{th:statistics}, \eqref{eq:M rc} and
\eqref{eq:X}.
\end{proof}

\section{Type $A_n^{(1)}$ dual bijection}
For this section we assume type $A_n^{(1)}$. We define and study the
properties of a dual analogue $\ddb$ of the $\db$ map that
corresponds to removing a tensor factor $B^{1\vee}$ from the left.
This is used to prove a duality symmetry (Theorem \ref{th:dualpath})
for the path-RC bijection in type $A_n^{(1)}$. This in turn is
useful for establishing the virtual bijections in section~\ref{sec:virtualbij}.

Let $\CCAD\subset \CCA$ be the category of tensor products of
crystals of the form $B^{1,s}$ and $B^{1,s\vee}$.

One goal of this section is to give a simpler way to compute $\phib$
for $B\in \CC^{A\vee}$. Since $\CCAD\subset\CCA$, Proposition
\ref{pp:phi} gives the definition of $\phib$. By \eqref{eq:dualKR}
$B^{1,s\vee}$ is isomorphic to $B^{n,s}$. The definition of $\phib$
involves left-splitting $B^{n,s}$, which produces columns $B^{n,1}$,
each of which have to be ``box split" into boxes $B^{1,1}$ and
removed by $\lh$.

We introduce a dual analogue $\ddb$ of $\db$, which removes an
entire column $B^{n,1}$ in a single step whose computation is
entirely similar to a single $\db$ (rather than $n$ of them).

Using $\ddb$, we can compute $\phib$ for $B\in \CCAD$ using
essentially single row techniques.

\subsection{Dual left hat}\label{ss:dual left hat}

Suppose that $B=B^{1\vee} \otimes B'$. In this particular case we
write $\ldh(B)=B'$. By Lemma \ref{lem:hwvhat} there is a map
$\ldh:P(B)\rightarrow P(\ldh(B))$ given by removing the left tensor
factor. Let $\Ldh$ be the multiplicity array of $\ldh(B)$.

The following algorithm is the same as $\db$ except that it starts from
large indices instead of small. The map $\ddb:RC(L)\rightarrow
RC(\Ldh)$ is defined as follows. Let $(\nu,J)\in RC(L)$. Initialize
$\ell^{(n+1)}=0$ and $\ell^{(0)}=\infty$. For $i$ from $n$ down to
$1$, assuming that $\ell^{(i+1)}$ has already been defined, let
$\ell^{(i)}$ be the smallest integer such that $(\nu,J)^{(i)}$ has a
singular string of length $\ell^{(i)}$ and $\ell^{(i)} \ge
\ell^{(i+1)}$. If no such singular string exists, let
$\ell^{(j)}=\infty$ for $1\le j\le i$. Let
$\rk^\vee(\nu,J)=(i+1)^\vee\in B^{1\vee}$ where $i$ is the maximum
index $i$ such that $\ell^{(i)}=\infty$.

\begin{example}
For $B=B^{1\vee}\otimes (B^1)^{\otimes 2} \otimes (B^2)^{\otimes 3} \otimes
(B^3)^{\otimes 2}$ of type $A_5^{(1)}$ and $\la=\Lab_1+\Lab_2+2\Lab_3+\Lab_4+\Lab_6$
the rigged configuration
\begin{equation*}
(\nu,J)= \quad \yngrc(3,0,2,0,2,0,2,0) \quad \yngrc(2,0,2,0,1,0) \quad
\yngrc(1,0,1,0) \quad \yngrc(1,1) \quad \yngrc(1,0)
\end{equation*}
is in $\RC(L,\la)$ with $L$ the multiplicity array corresponding to $B$.
The same configuration now written with the vacancy number next to each part is
\begin{equation*}
\qquad \yngrc(3,1,2,1,2,1,2,1) \quad \yngrc(2,0,2,0,1,0) \quad
\yngrc(1,0,1,0) \quad \yngrc(1,1) \quad \yngrc(1,0)
\end{equation*}
Then
\begin{equation*}
\ddb(\nu,J)=\quad \yngrc(3,0,2,0,2,0,2,0) \quad \yngrc(2,0,2,0) \quad \yngrc(1,0)
\quad \es \quad \es
\end{equation*}
and $\rk^\vee(\nu,J)=2^\vee$.
\end{example}

Given $\mu\in\lm$, there is also an inverse of the dual algorithm $\ddb$
associated with the weight $(\la-\mu)^\vee$ similar to the inverse of
$\db$ as defined in section~\ref{ss:deltas}.

\begin{prop} \label{pp:dbv} $\ddb:RC(L)\rightarrow RC(\Ldh)$ is a
well-defined injective map such that the diagram commutes:
\begin{equation} \label{eq:dualhatdiag}
\begin{CD}
  P(B) @>{\phib}>> RC(L) \\
  @V{\ldh}VV @VV{\ddb}V \\
  P(\ldh(B)) @>>{\phib}> RC(\Ldh).
\end{CD}
\end{equation}
Moreover, if $\phib(b_1\otimes b)=(\nu,J)$ then
 $\phib(b)=\ddb(\nu,J)$ and $b_1=\rk^\vee(\nu,J)$.
\end{prop}

\begin{proof} The map $\ldh$ removes $B^{1\vee}\cong B^{n,1}$.
This may be achieved by $n$ applications of $\lh\circ \llb$, which
splits a box from a column and then removes it. Let $\Delta$ be the
corresponding $n$-fold composition of maps $\db \circ \rclb$. It
must be shown that $\Delta(\nu,J)=\ddb(\nu,J)$.

Let $a^\vee=b_1$. The letters $1\le a_1<a_2<\cdots<a_n\le n+1$ in
$b_1$ (where $b_1$ is viewed as an element of $B^{n,1}$) satisfy
\begin{equation*}
a_i=\begin{cases} i & \text{for $1\le i<a$}\\
i+1 & \text{for $a\le i\le n$.} \end{cases}
\end{equation*}
It is clear from Proposition~\ref{pp:rcbs} and the algorithm for
$\db$ of section~\ref{ss:deltas} that the composition $\db\circ
\rclb$ corresponding to the letter $a_i$ for $a\le i\le n$ in $b_1$
removes a box from one string of length $s^{(i)}$ in the $i$-th
rigged partition and leaves all other strings unchanged. It also
follows from the algorithms and change of vacancy numbers that
$s^{(i)}\ge s^{(i+1)}$ and that $s^{(i)}$ is the length of the
smallest singular string in $(\nu,J)^{(i)}$ with this property. For
$1\le i<a$ the composition $\db\circ \rclb$ leaves the rigged
configuration unchanged with $s^{(i)}=\infty$. It follows by
induction that $s^{(i)}=\ell^{(i)}$, where $\ell^{(i)}$ as in the
definition of $\ddb$, and hence that $\Delta(\nu,J)=\ddb(\nu,J)$.
\end{proof}

\subsection{Dual left split} We restate left splitting for a
special case. Suppose $B=B^{s\vee} \otimes B'$ for $s\ge 2$. Define
$\ls^\vee:B^{s\vee}\rightarrow \ls^\vee(B):=B^{1\vee}\otimes
B^{s-1\vee}$ to be the composite map
\begin{equation*}
\begin{CD}
B^{s\vee} @>{\sim}>> B^{n,s} @>{\ls}>> B^{n,1} \otimes B^{n,s-1}
@>{\sim}>> B^{1\vee} \otimes B^{s-1\vee}.
\end{CD}
\end{equation*}
By Example \ref{ex:dualArow} we may write $b\in B^{s\vee}$ as a word
of length $s$ in the dual alphabet. Computing $\ls^\vee$ using
Example \ref{ex:dualA}, it is seen that $\ls(b)=b_2\otimes b_1$
where $b_2$ is the leftmost dual letter in $b$ and $b_1$ is the
remaining word of length $s-1$ in the dual alphabet.

Let $\Lds$ be the multiplicity array for $\ls^\vee(B)$. Let us
denote by $\rcls^\vee$ the map on RCs which corresponds to
$\ls^\vee$ under the path-RC bijection $\phib$. It is the map
$\rcls$ with respect to $B^{n,s}$ and is therefore inclusion (with
some changes in vacancy numbers). With these definitions the
following diagram commutes by Proposition \ref{pp:phi} for
$A_n^{(1)}$:
\begin{equation} \label{eq:dualsplitdiag}
\begin{CD}
P(B,\la) @>{\phib}>> \RC(L,\la) \\
@V{\ls^\vee}VV @VV{\rcls^\vee}V \\
P(\ls^\vee(B),\la) @>>{\phib}> \RC(\Lds,\la).
\end{CD}
\end{equation}

\subsection{The bijection $\phib$ for $\CCAD$}
The results of this section to this point may be summarized as
follows.

\begin{prop}\label{pp:singledualbij}
There is a unique bijection $\phib:P(B)\rightarrow \RC(L)$
satisfying the following properties. It sends the empty path to the
empty rigged configuration, and if the leftmost tensor factor in $B$
is:
\begin{enumerate}
\item $B^1$: \eqref{eq:hatdiag} holds.
\item $B^s$ for $s\ge2$: \eqref{eq:splitdiag} holds.
\item $B^{1\vee}$: \eqref{eq:dualhatdiag} holds.
\item $B^{s\vee}$ for $s\ge2$: \eqref{eq:dualsplitdiag} holds.
\end{enumerate}
\end{prop}

\subsection{Duality on paths and the bijection}
Let $B\in \CCAD$.

Let $L$ and $L^\vee$ be the multiplicity arrays for $B$ and $B^\vee$
respectively. Explicitly, $L^{\vee(a)}_i = L^{(n+1-a)}_i$ for $1\le
a \le n$. Given a classical highest weight $\la$, let
$\la^\vee=-w_0\la$ be the highest weight of the contragredient dual
module to the $A_n$-module highest weight $\la$. There is a
bijection $\vee:\RC(L,\la)\rightarrow \RC(L^\vee,\la^\vee)$ given by
$(\nu,J)\mapsto (\nu',J')$ where ${\nu'}^{(a)}=\nu^{(n+1-a)}$ and
${J'}^{(a,i)}$ is obtained from $J^{(n+1-a,i)}$ by complementation
within the $m_i^{(n+1-a)}(\nu) \times p_i^{(n+1-a)}(\nu)$ rectangle.

\begin{theorem} \label{th:dualpath} \cite{OSS:2003a} Let $B\in\CCA$,
$B^\vee$ its contragredient dual, and $L$ and $L^\vee$ their
respective multiplicity arrays. The diagram commutes:
\begin{equation*}
\begin{CD}
  P(B) @>{\phib}>> RC(L) \\
  @V{\vee}VV @VV{\vee}V \\
  P(B^\vee) @>>{\phib}> RC(L^\vee).
\end{CD}
\end{equation*}
\end{theorem}

\section{Virtual bijection}
\label{sec:virtualbij}

In this section we will prove $X=M$ for the category $\CC$ for the
nonsimply-laced algebras. For the simply-laced types $A_n^{(1)}$ and
$D_n^{(1)}$ this was proved in Corollary~\ref{cor:X=M AD}. For the
non-simply-laced affine families it suffices to prove the following
theorem.

\begin{theorem} \label{th:virtualbij} For $B\in \CC$, let
$\Psi:B\rightarrow \Vh$ be the virtual crystal embedding, $L$ and
$\widehat{L}$ the multiplicity arrays for $B$ and $\Vh$
respectively. Then the simply-laced bijection
$\phib_{\Lhat}:P(\Vh)\rightarrow \RC(\Lhat)$ restricts to a
bijection $\vphib:\Pv(B)\rightarrow \RCv(L)$.
\end{theorem}

As an immediate corollary we obtain:
\begin{corollary}\label{cor:X=VX=VM=M}
For $\la\in\Pb^+$, $B\in\CC$ and $L$ the corresponding multiplicity array
we have
\begin{equation*}
X_{B,\la}(q)=VX_{B,\la}(q)=VM_{L,\la}(q)=M_{L,\la}(q).
\end{equation*}
\end{corollary}
\begin{proof}
The left and right equalities were proven in Theorem~\ref{th:X=VX} and
Corollary~\ref{cor:M=VM}, respectively. The middle equality follows from
Theorems~\ref{th:statistics} and~\ref{th:virtualbij}.
\end{proof}
The remainder of this section is occupied with the proof of
Theorem~\ref{th:virtualbij}.

\subsection{Virtual $\lh$}\label{sec:virtual lh}
Suppose $B=B_X=B^1_X \otimes B'_X\in \CC$ with virtual crystal
embeddings $\emb:B_X\rightarrow \Vh$ and $\emb':B'_X\rightarrow
\Vh'$. By abuse of notation we write $\vlh(\Vh)=\Vh'$. The map
$\vlh:\Vh\rightarrow\Vh'$ is defined by
\begin{enumerate}
\item If $Y=A_{2n-1}^{(1)}$ then $\vlh:B^{1\vee}_Y \otimes B^1_Y
\otimes \Vh'\rightarrow\Vh'$ is defined by $\vlh = \lh \circ \ldh$,
which drops the two leftmost factors in $\Vh$.
\item If $Y=D_{n+1}^{(1)}$ and $X=B_n^{(1)}$ then $\vlh:B_Y^2 \otimes
\Vh'\rightarrow \Vh'$ is defined by $\vlh= \lh \circ \lh \circ \ls$.
This accomplishes the same thing as deleting the tensor factor
$B_Y^2$.
\item If $Y=D_{n+1}^{(1)}$ and $X=A_{2n-1}^{(2)}$ then $\vlh:B_Y^1
\otimes \Vh'\rightarrow\Vh'$ is defined by $\vlh=\lh$.
\end{enumerate}
Note that in each case the total effect of the map $\vlh:\Vh^1
\otimes \Vh'\rightarrow \Vh'$ is to drop the tensor factor $\Vh^1$.
Therefore the following diagram commutes trivially:
\begin{equation*}
\begin{CD}
B_X^1 \otimes B_X' @>{\emb\otimes\emb}>> \Vh^1 \otimes \Vh' \\
@V{\lh}VV @VV{\vlh}V \\
B_X' @>>{\emb}> \Vh'.
\end{CD}
\end{equation*}

\subsection{Virtual $\ls$} \label{sec:virtual ls}
Let $s\ge2$. Recall the virtual $\rs$ map $\vrs:\Vh^s\rightarrow
\Vh^{s-1} \otimes \Vh^1$ defined in the proof of Proposition
\ref{pp:vrs}. Define the virtual $\ls$ map $\vls:\Vh^s \rightarrow
\Vh^1 \otimes \Vh^{s-1}$ by
\begin{equation} \label{eq:vlsdef}
\vls=*\circ \vrs \circ *.
\end{equation}

\begin{prop} \label{pp:xvls} The map $\vls:\Vh^s\rightarrow \Vh^1\otimes \Vh^{s-1}$
is described explicitly as follows.
\begin{enumerate}
\item If $Y=A_{2n-1}^{(1)}$ then $\vls:B_Y^{s\vee}\otimes B_Y^s\rightarrow
B_Y^{1\vee} \otimes B_Y^1 \otimes B_Y^{s-1\vee}\otimes B_Y^{s-1}$ is
the composition
\begin{equation*}
\begin{split}
 &B_Y^{s\vee} \otimes B_Y^s \stackrel{\ls_Y^\vee\otimes 1}{\longrightarrow}
B_Y^{1\vee} \otimes B_Y^{s-1\vee} \otimes
B_Y^s\stackrel{R}{\longrightarrow}
B_Y^s \otimes B_Y^{1\vee} \otimes B_Y^{s-1\vee}\\
\stackrel{\ls_Y\otimes 1\otimes 1}{\longrightarrow}& B_Y^1 \otimes
B_Y^{s-1} \otimes B_Y^{1\vee} \otimes B_Y^{s-1\vee}
\stackrel{R}{\longrightarrow} B_Y^{1\vee}\otimes B_Y^1 \otimes
B_Y^{s-1\vee} \otimes B_Y^{s-1}.
\end{split}
\end{equation*}
\item If $Y=D_{n+1}^{(1)}$ and $X=B_n^{(1)}$ then
$\vls:B_Y^{2s}\rightarrow B_Y^2\otimes B_Y^{2s-2}$ is the map that
splits off the first two symbols, that is, $uv\mapsto u\otimes v$
where $uv\in B_Y^{2s}$, $u\in B_Y^2$, and $v\in B_Y^{2s-2}$.
\item If $Y=D_{n+1}^{(1)}$ and $X=A_{2n-1}^{(2)}$ then define
$\vls=\ls_Y:B_Y^s\rightarrow B_Y^1\otimes B_Y^{s-1}$.
\end{enumerate}
\end{prop}
\begin{proof} It is enough to check these on
classical highest weight vectors. This is easy because the various
crystals are multiplicity-free as classical crystals.
\end{proof}

\begin{remark} \label{rem:vls}
Let $B=B_X=B^s \otimes B'$ and let $\Psi:B\rightarrow \Vh$ and
$\Psi':B'\rightarrow \Vh'$ be the virtual crystal realizations. By
abuse of notation we write $\vls(\Vh)=\Vh^1 \otimes \Vh^{s-1}
\otimes \Vh'$. We also use the notation $\vls$ for the map $\vls
\otimes 1_{\Vh'}:\Vh^s\otimes \Vh'\rightarrow \Vh^1\otimes\Vh^{s-1}
\otimes \Vh'$. It also satisfies \eqref{eq:vlsdef}.
\end{remark}

\subsection{Virtual $\db$ and $\rcls$} Given the virtual crystal
embedding $\emb:B_X\rightarrow\Vh$, let $L$ and $\Lhat$ be the
multiplicity arrays for $B_X$ and $\Vh$ respectively. The maps
$\vdb$ and $\vrcls$ are defined to be the maps on rigged
configurations which correspond under $\phib$ to the maps $\vlh$ and
$\vls$. More precisely, since $\phib$ is a bijection for type $Y$
there are unique maps $\vdb$ and $\vrcls$ defined by the commutation
of the diagrams
\begin{equation}\label{eq:vdelta vj}
\begin{CD}
  P(\Vh) @>{\phib_{\Lhat}}>> RC(\Lhat) \\
  @V{\vlh}VV @VV{\vdb}V \\
  P(\vlh(\Vh)) @>>{\phib_{\vlh(\Lhat)}}> RC(\vlh(\Lhat))
\end{CD}
\qquad\text{and}\qquad
\begin{CD}
  P(\Vh) @>{\phib_{\Lhat}}>> RC(\Lhat) \\
  @V{\vls}VV @VV{\vrcls}V \\
  P(\vls(\Vh)) @>>{\phib_{\vls(\Lhat)}}> RC(\vls(\Lhat))
\end{CD}
\end{equation}
where $\vlh(\Lhat)$ and $\vls(\Lhat)$ are the multiplicity arrays for
$\vlh(\Vh)$ and $\vls(\Vh)$ respectively.

For $\la\in\Pb^+(X)$ let $\rkh:\RC(\Lhat,\emb(\la))\rightarrow
\Vh^1$ be the map which gives the tensor product of the ranks of the
sequence of rigged configurations that occur during the computation
of $\vdb$.

\begin{lemma} \label{lem:vdb} $\vdb$ maps $\RCv(L)$ into
$\RCv(\Lh)$ and $\rkh$ maps $\RCv(L)$ into
$\Image(\emb:B_X^1\rightarrow \Vh^1)$.
\end{lemma}
\begin{proof} The proof proceeds by cases.
\subsubsection*{$X=C_n^{(1)}$ and $Y=A_{2n-1}^{(1)}$.} According to
Definition~\ref{def:VRC} the elements $(\nh,\Jh)\in\RCv(L)$ have the
following properties:
\begin{enumerate}
\item \label{sym}$\mh_i^{(a)}=\mh_i^{(2n-a)}$ and
$\Jh^{(a,i)}=\Jh^{(2n-a,i)}$;
\item \label{m even}$\mh^{(n)}_i=0$ if $i$ is odd;
\item \label{J even} The parts of $\Jh^{(n,i)}$ are even.
\end{enumerate}
{}From \eqref{eq:vdelta vj} and Proposition~\ref{pp:dbv} it is clear
that $\vdb=\db \circ \ddb$. It must be shown that $\vdb(\nh,\Jh)$
also possesses the three properties \eqref{sym}-\eqref{J even}. Let
$\ell^{\vee (a)}$ the lengths of the strings selected by $\ddb$ and
$\ell^{(a)}$ be the lengths of the strings selected by the
subsequent application of $\db$. Let
$\rk^\vee(\nh,\Jh)=(2n+1-r)^\vee$ for some $1\le r\le 2n$. If $r\le
n$, it is clear from the definitions that $\ell^{(a)}=\ell^{\vee
(2n-a)}$ for $1\le a<r$, so that points \eqref{sym}-\eqref{J even}
still hold. Here $\rkh(\nh,\Jh)=(2n+1-r)^\vee \otimes r=\emb(r)$.
For $r=n+1$, we must have $\ell^{\vee(n+1)}<\ell^{\vee(n)}$ since
otherwise by the symmetry \eqref{sym}
$\ell^{\vee(n-1)}=\ell^{\vee(n)}=\ell^{\vee(n+1)}<\infty$ which
contradicts the assumption that $r=n+1$. However, this implies that
$\ell^{(a)}=\ell^{\vee (2n-a)}$ for $1\le a<n$ and $\ell^{
(n)}=\ell^{\vee(n)}-1$. Since the vacancy numbers are all even
\eqref{sym}-\eqref{J even} remain valid. One has
$\rkh(\nh,\Jh)=n^\vee \otimes (n+1)=\emb(\overline{n})$. Finally let
$r>n+1$ and let $r'\le n$ be minimal such that
$\ell^{\vee(2n-r')}=\ell^{\vee(n)}$. By symmetry \eqref{sym} we have
$\ell^{\vee(a)}=\ell^{\vee(n)}$ for all $r'\le a\le 2n-r'$. By the
algorithms for $\ddb$ and $\db$ and properties \eqref{sym}-\eqref{J
even} for $(\nh,\Jh)$ it follows that $\ell^{(a)}=\ell^{\vee
(2n-a)}$ for $1\le a<r'$ and $2n-r'<a<r$, and $\ell^{(a)}=\ell^{\vee
\
(2n-a)}-1$ for $r'\le a\le 2n-r'$. Again this implies that
properties \eqref{sym}-\eqref{J even} hold for $\vdb(\nh,\Jh)$. Then
$\rkh(\nh,\Jh)=(2n+1-r)^\vee \otimes r=\emb(\overline{2n+1-r})$.

\subsubsection*{$X=A_{2n}^{(2)}$ and $Y=A_{2n-1}^{(1)}$.} The elements in
$\RCv(L)$ are characterized by points \eqref{sym} and \eqref{J
even}. Everything goes through as for the case $X=C_n^{(1)}$ except
that, since $\nh^{(n)}$ may contain odd parts, it is possible that
$\ell^{\vee(n)}=1$. In this case $\ell^{\vee(a)}=1$ for all $1\le
a\le 2n-1$ by point \eqref{sym}. Then $\ell^{(a)}=\infty$ for all
$1\le a\le 2n-1$, so that $\vdb(\nh,\Jh)$ again satisfies
\eqref{sym} and \eqref{J even}. Then $\rkh(\nh,\Jh)=1^\vee\otimes
1=\emb(\vn)$.

\subsubsection*{$X=D_{n+1}^{(2)}$ and $Y=A_{2n-1}^{(1)}$.} The elements in
$\RCv(L)$ are characterized by point \eqref{sym}. The proof goes
through as before except that $\Jh^{(n,i)}$ could have an odd part.
This could only change the computation of $\db\circ\ddb$ if such an
odd part were selected. Recall that $p_i^{(n)}$ is even for all $i$.
Therefore the odd part cannot be selected by $\ddb$. It can only be
selected by $\db$ if $\rk^\vee(\nh,\Jh)=(n+1)^\vee$ and the odd part
has size $p_i^{(n)}-1$ for some $i\ge \ell^{\vee(n+1)}$. By point
\eqref{sym} and the fact that $(\nh,\Jh)^{(a)}$ is unchanged by
$\ddb$ for $1\le a\le n-1$, we have $\ell^{(a)}=\ell^{\vee(2n-a)}$
for $1\le a \le n-1$ and $\ell^{(n)}$ is the odd (now singular)
part. Thus after applying $\db\circ\ddb$ point \eqref{sym} still
holds. $\rkh(\nh,\Jh)=(n+1)^\vee \otimes (n+1)=\emb(0)$. Note that
$\ell^{(n+1)}=\infty$ since $\ddb$ caused the strings in the
$(n+1)$-th rigged partition that were longer than
$\ell^{\vee(n+1)}$, to become nonsingular.

\subsubsection*{$X=A_{2n}^{(2)\dagger}$ and $Y=A_{2n-1}^{(1)}$.}
The elements in $\RCv(L)$ are characterized by \eqref{sym} and
\begin{enumerate}
\item[(3')] The parts of $\Jh^{(n,i)}$ have the same parity as $i$.
\end{enumerate}
Let $\rk^\vee(\nh,\Jh)=(2n+1-r)^\vee$ for some $1\le r\le 2n$. If
$r\le n$, we have as for the case $X=C_n^{(1)}$ that
$\ell^{(a)}=\ell^{\vee (2n-a)}$ for $1\le a<r$, so that \eqref{sym},
and (3') still hold and $\rkh(\nh,\Jh)=(2n+1-r)^\vee \otimes
r=\emb(r)$.

If $r=n+1$, note that $\ell^{\vee(n)}\in 2\Z$ since all vacancy
numbers $p_i^{(n)}$ are even, so that by (3') only the riggings for
$i$ even can possibly be singular. As in case $C_n^{(1)}$ we must
have $\ell^{\vee(n+1)}<\ell^{\vee(n)}$. By symmetry \eqref{sym} we
have $\ell^{(a)}=\ell^{\vee(2n-a)}$ for $1\le a<n$. The application
of $\ddb$ changes the vacancy numbers in the $n$-th rigged partition
corresponding to the strings of length $i$ for $\ell^{\vee(n+1)}\le
i<\ell^{\vee(n)}$ by $-1$, which makes these vacancy numbers odd. In
particular, the rigging of the new string of length
$\ell^{\vee(n)}-1$ is odd. In addition, $\ell^{\vee(n+1)}\le
\ell^{(n)}<\ell^{\vee(n)}$ and by (3') $\ell^{(n)}$ must be odd. By
the change in vacancy number after the application of $\db$, the new
rigging of the string of length $\ell^{(n)}-1$ must be even. Hence
\eqref{sym} and (3') hold for $\vdb(\nh,\Jh)$ and
$\rkh(\nh,\Jh)=(n+1)^\vee \otimes (n+1)=\emb(0)$.

If $r>n+1$, let $r'\le n$ be defined as for the case $C_n^{(1)}$. As
before $\ell^{\vee(a)}=\ell^{\vee(n)}$ for $r'\le a\le 2n-r'$. If
$r'<n$ everything goes through as in case $C_n^{(1)}$. If $r'=n$
(which means that $\ell^{\vee(n+1)}<\ell^{\vee(n)}$), by the same
arguments as for $r=n+1$, we have $\ell^{(a)}=\ell^{\vee(2n-a)}$ for
$a\neq n$, $\ell^{\vee(n+1)}\le \ell^{(n)}<\ell^{\vee(n)}$ and (3')
holds for the new riggings. Hence properties \eqref{sym} and (3')
hold for $\vdb(\nh,\Jh)$ and $\rkh(\nh,\Jh)=(2n+1-r)^\vee \otimes
r=\emb(\overline{2n+1-r})$.

\subsubsection*{$X=B_n^{(1)}$ and $Y=D_{n+1}^{(1)}$.} The elements in $\RCv(L)$ are
characterized by
\begin{enumerate}
\item \label{sym D}$m_i^{(n)}=m_i^{(n+1)}$ and $J^{(n,i)}=J^{(n+1,i)}$ for all $i>0$;
\item \label{even D}$\nu^{(a)}$ and $J^{(a,i)}$ have only even parts for $1\le a<n$.
\end{enumerate}
By section~\ref{sec:virtual lh} and \eqref{eq:vdelta vj} we have
$\vdb=\db\circ\db\circ\rcls$. Let $\ell^{(a)}$ and $\lb^{(a)}$
(resp. $s^{(a)}$ and $\sbar^{(a)}$) be the length of the selected
strings for the right (resp. left) $\db$. Then it follows from the
definition of $\rcls$, $\db$ and point \eqref{even D} that
$s^{(a)}=\ell^{(a)}-1$ for $1\le a<n$. Furthermore from point \eqref{sym D} we
obtain that $\ell^{(n)}=\ell^{(n+1)}>s^{(n)}=s^{(n+1)}$, and again
by point \eqref{even D} that $\sbar^{(a)}=\lb^{(a)}-1$ for $1\le a<n$. This implies
that points \eqref{sym D} and \eqref{even D} hold for
$\vdb(\nh,\Jh)$. Moreover, let $x=\rk(\nh,\Jh)$ and $y=\rk(\db(\nh,\Jh))$.
Note that $x,y\neq n+1,\overline{n+1}$ because of point (1).
Also $x=y$ except possibly $x=n$ and $y=\overline{n}$. Then
$\rkh(\nh,\Jh)=xx=\emb(x)$ if $x=y$ or $\rkh(\nh,\Jh)=n\overline{n}=\emb(0)$
if $x\neq y$.

\subsubsection*{$X=A_{2n-1}^{(2)}$ and $Y=D_{n+1}^{(1)}$.} The elements
in $\RCv(L)$ are characterized by point \eqref{sym D}. It is obvious
from its definition that $\vdb=\db$ preserves this property. Let
$x=\rk(\nh,\Jh)$. As before $x\neq n+1,\overline{n+1}$ because of point (1).
Then $\rkh(\nh,\Jh)=x=\emb(x)$.
\end{proof}

Thus we may define the virtual rank map $\vrk:\RCv(L)\rightarrow
B_X^1$ by $\vrk(\nh,\Jh)=x$ where $\emb(x)=\rkh(\nh,\Jh)$ for all
$(\nh,\Jh)\in\RCv(L)$. Then we have:

\begin{prop} \label{pp:vboxbij} The map
$(\vdb,\vrk):\RCv(L,\la)\rightarrow \bigcup_{\mu\in\lm}
\RCv(\Lh,\mu) \times B_X^1$ is injective.
\end{prop}

For the proof of Theorem \ref{th:virtualbij} we also need the inverse
to Lemma~\ref{lem:vdb} which involves the inverse of $\vdb$.
Let $\la\in \Pb^+_X$, $L=(L_1,L_2,\ldots)$ a multiplicity array and
$\Lp=(L_1+1,L_2,L_3,\ldots)$. Denote by $\RCtv(L,\la)$ the subset of
$\RCv(L,\la)\times B^1$ given by $((\nu,J),b)$ such that $\la+\wt(b)\in\Pb^+$
and if $b=0$ then also $\la_n>0$. Let $\hat{b}=\emb(b)$. By abuse of
notation we define
\begin{equation*}
\vdb^{-1}:\RCtv(L,\la)\to\bigcup_{\beta\in\lp}\RC(\widehat{\Lp},\emb(\beta)).
\end{equation*}
If $Y=A_{2n-1}^{(1)}$, let $\hat{b}=b_1\otimes b_2$.
Then $\vdb^{-1}((\nu,J),b)={\ddb}^{-1}(\db^{-1}((\nu,J),b_2),b_1)$, with
$\db^{-1}$ as defined in section~\ref{ss:deltas} and ${\ddb}^{-1}$ as defined in
section~\ref{ss:dual left hat}.
If $Y=D_{n+1}^{(1)}$ and $X=B_n^{(1)}$, let $\hat{b}=xy$. Then
$\vdb^{-1}((\nu,J),b)={\db}^{-1}(\db^{-1}((\nu,J),y),x)$.
Finally for $Y=D_{n+1}^{(1)}$ and $X=A_{2n-1}^{(2)}$, let $\hat{b}=x$. Then
$\vdb^{-1}((\nu,J),b)=\db^{-1}((\nu,J),x)$.

\begin{lemma}\label{lem:vdb inv}
Given $\la$, $L$, $\Lp$, $b$ and $\hat{b}$ as above
the map $\vdb^{-1}$ maps $\RCtv(L,\la)$ into $\bigcup_{\beta\in\lp}\RCv(\Lp,\beta)$.
\end{lemma}
\begin{proof}
The proof is very similar to the proof of Lemma~\ref{lem:vdb}.
\end{proof}

\begin{lemma} \label{lem:vls} $\vrcls$ maps $\RCv(L)$ into
$\RCv(\Ls)$.
\end{lemma}
\begin{proof} Let $Y=A_{2n-1}^{(1)}$. By \eqref{eq:vdelta vj} and
section~\ref{sec:virtual ls} we have $\vrcls=\rcls\circ \rclds$.
Both $\rcls$ and $\rclds$ are inclusions that do not change the
rigged configuration (only certain vacancy numbers). Hence if
$(\nh,\Jh)\in\RCv(L)$ has the characterization as stated in the
previous lemma, then so does $\vrcls(\nh,\Jh)$.

Let $Y=D_{n+1}^{(1)}$. If $X=A_{2n-1}^{(2)}$, we have
$\vrcls=\rcls$. For $X=B_n^{(1)}$, let $B=B^s\otimes B'$ for $s\ge2$
and embeddings $\Psi:B\rightarrow\Vh$ and $\Psi':B'\rightarrow\Vh'$
with $\Vh=B_Y^{2s} \otimes \Vh'$. It can be shown (using $*$ and
properties of $\rs$) that if $x,y\in B_Y^1$ and $u\in B_Y^{2s-2}$
are such that $xyu\in B_Y^{2s}$ then for any $b'\in \Vh'$ one has
$\vls(xyu\otimes b')=xy \otimes u \otimes b'$. One may show that the
corresponding operation on RCs is inclusion. This may be seen by
observing that $\ls \circ \vls:B_Y^{2s}\otimes \Vh' \rightarrow
B_Y^1 \otimes B_Y^1\otimes B_Y^{2s-2} \otimes \Vh'$, which sends
$xyu\otimes b'$ to $x\otimes y\otimes u\otimes b'$, can also be
computed by a composition of $\ls$ maps and $R$-matrices, whose
corresponding maps on RCs are inclusions. This proves that
$\vrcls(\nh,\Jh)\in\RCv(\Ls)$.
\end{proof}

\subsection{Proof of Theorem \ref{th:virtualbij}}
It must be shown that the bijection
$\phib_{\Lhat}:P(\Vh)\rightarrow \RC(\Lhat)$ maps $\Pv(B)$ (1)
into and (2) onto $\RCv(L)$, thereby defining a bijection
$\vphib_L:\Pv(B)\rightarrow\RCv(L)$ by restriction. Let $B=B^s
\otimes B'$ with $\emb:B\rightarrow \Vh^s\otimes \Vh'$.

\subsubsection*{The case $s=1$:}
For (1) consider a typical element of $\Pv(B,\la)$, given by
$\emb(b)$ with $b\in P(B,\la)$. Write $b=x\otimes b'$ with $x\in
B_X^1$ and $b'\in P(B',\mu)$. Then $\emb(b')\in\Pv(\lh(B),\mu)$. Let
$(\nh,\Jh)=\phib_{\Lhat}(\emb(b))\in\RC(\Lhat)$. It must be shown
that $(\nh,\Jh)\in \RCv(L,\la)$. By \eqref{eq:vdelta vj} and
induction one has $\vdb(\nh,\Jh)\in \RCv(\Lh,\mu)$ and
$\rkh(\nh,\Jh)=\emb(x)$. By Lemma~\ref{lem:vdb inv} we can conclude
that $(\nh,\Jh)\in\RCv(L,\la)$.

For (2) let $(\nh,\Jh)\in\RCv(L)$. Let $\bh=\xh\otimes \bhp\in
P(\Vh)$ (with $\xh\in\Vh^1$ and $\bhp\in \Vh'$) be such that
$\phib_{\Lhat}(\bh)=(\nh,\Jh)$. It must be shown that $\bh\in
\Pv(B)$. By \eqref{eq:vdelta vj} we have
$\phib_{\vlh(\Lhat)}(\vlh(\bh))= \vdb(\phib_{\Lhat}(\bh))=
\vdb(\nh,\Jh)\in\RCv(\vlh(\Lhat))$. By induction $\bhp=\vlh(\bh)\in
\Pv(\lh(B))$; write $\bhp=\emb(b')$ for some $b'\in B'$. By Lemma
\ref{lem:vdb} and \eqref{eq:vdelta vj}, $\xh=\emb(x)$ for
$x=\vrk(\nh,\Jh)$. Let $b=x \otimes b'\in B$. By definition
$\emb(b)=\emb(x)\otimes \emb(b')=\xh\otimes \bhp=\bh$. Therefore
$\bh\in \Pv(B)$ as desired.

\subsubsection*{The case $s\ge2$:}
For (1), a typical element of $\Pv(B)$ has the form $\emb(b)$ for
$b\in P(B)$. Let $\phib_{\Lhat}(\emb(b))=(\nh,\Jh)\in \RC(\Lhat)$.
It must be shown that $(\nh,\Jh)\in \RCv(L)$. Note that
$\vrcls(\nh,\Jh)=\vrcls(\phib_{\Lhat}(\emb(b)))=
\phib_{\Lhat^s}(\vls(\emb(b)))\in\RCv(\Ls)$
by \eqref{eq:vdelta vj} and induction. But
$\vrcls(\nh,\Jh)=(\nh,\Jh)$ and $(\nh,\Jh)\in\RC(\Lhat)$. It follows
that $(\nh,\Jh)\in\RCv(L)$.

For (2), let $(\nh,\Jh)\in\RCv(L)$. Let $\bh\in P(\Vh)$ be such that
$\phib_{\Lhat}(\bh)=(\nh,\Jh)$. It must be shown that $\bh\in
\Pv(B)$. By \eqref{eq:vdelta vj} and induction,
$\vrcls(\nh,\Jh)=\vrcls(\phib_{\Lhat}(\bh))=\phib_{\Lhat^s}(\vls(\bh))\in
\RCv(\Ls)$. Therefore $\vls(\bh)\in \Pv(\ls(B))$. We conclude that
$\bh\in \Pv(B)$ by \eqref{eq:vlsdef}, Proposition \ref{pp:vrs} point
\ref{it:vres}, and Proposition \ref{pp:*emb}.

This concludes the proof of Theorem \ref{th:virtualbij}.

\end{document}